\documentclass[12pt,epsfig,amsfonts]{amsart}
   \newtheorem{lemma}{Lemma}[section]
   \newtheorem{theorem}{Theorem}
   
   \newtheorem{prop}{Proposition}[section]

\usepackage{epsfig}

\numberwithin{equation}{section}

\begin{document}

\pagestyle{plain}
\author{Qiudong Wang}
\address{Department of Mathematics\\
The University of Arizona\\
Tucson, Arizona 85721-0089} \email[Qiudong
Wang]{dwang@math.arizona.edu} \urladdr[Qiudong
Wang]{math.arizona.edu/\textasciitilde dwang}

\date{January, 2008}

\begin{abstract}
In this paper we use the second order equation
\begin{equation*}
\frac{d^2 q}{dt^2} + (\lambda - \gamma q^2) \frac{d q}{dt} - q +
q^3 = \mu q^2 \sin \omega t
\end{equation*}
as a demonstrative example to illustrate how to apply the analysis
of \cite{WO} and \cite{WOk} to the studies of concrete equations.
We prove, among many other things, that there are positive measure
sets of parameters $(\lambda, \gamma, \mu, \omega)$ corresponding
to the case of intersected (See Fig. 1(a)) and the case of
separated (See Fig. 1(b)) stable and unstable manifold of the
solution $q(t) = 0$, $t \in \mathbb R$ respectively, so that the
corresponding equations admit strange attractors with SRB
measures.
\end{abstract}

\title{On the dynamics of a time-periodic equation}

\maketitle

In the history of the theory of dynamical systems, ordinary
differential equations have served as a source of inspirations and
a test ground. There are mainly two ways in relating the studies
of differential equations to the studies of maps, both originated
from the work of H. Poincar\'e. The first is to use return maps
locally defined on Poincar\'e sections for the studies of
autonomous equations. The second is to use the globally defined
time-T maps for the studies of time-periodic equations.

The studies of periodically forced second order equations, such as
van de Pol's equation, Duffing's equation and the equation for
non-linear pendulums, have played substantial roles in shaping the
chaos theory in modern times (\cite{Ar}, \cite{D}, \cite{V},
\cite{Lev1}, \cite{Le}, \cite{GH}). When a homoclinic solution is
periodically perturbed, the stable and the unstable manifold of
the perturbed saddle intersect each other within a certain range
of forcing parameters, generating homoclinic tangles and chaotic
dynamics (\cite{P}, \cite{S}, \cite{M}). See Fig. 1(a). There are
also forcing parameters for which the stable and the unstable
manifold of the perturbed saddle are pulled apart. See Fig. 1(b).

\begin{picture}(12, 4)
\put(2,0){ \psfig{figure=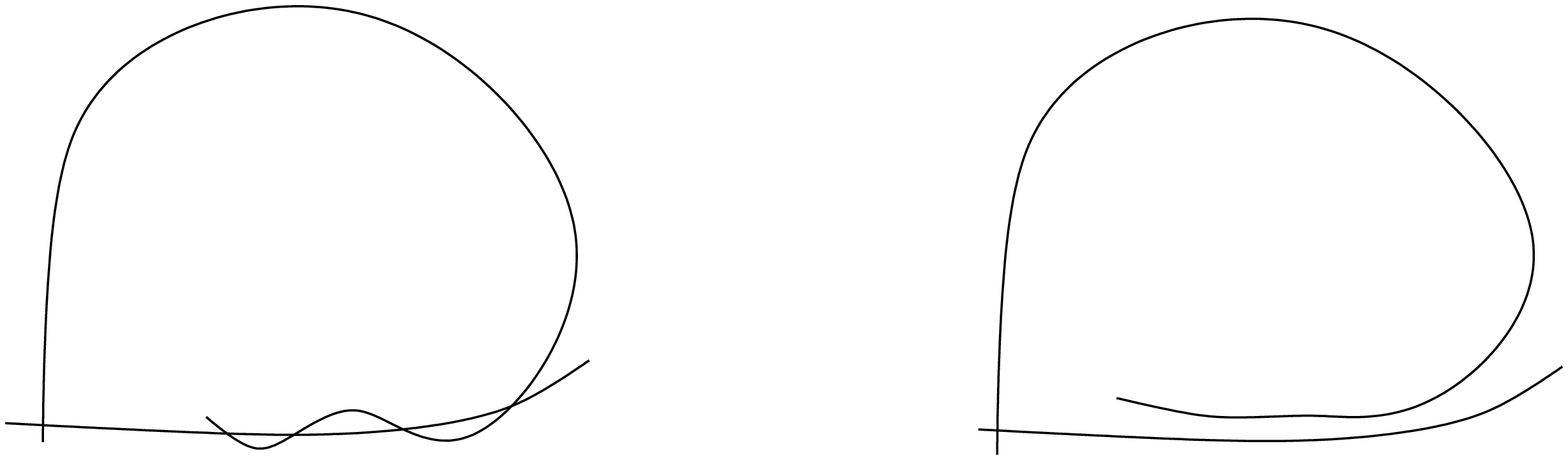,height = 3cm, width =
11cm} }
\end{picture}

\centerline{\ \ \ (a)  \hspace{2.2in} (b)}

\smallskip

\centerline{Fig. 1 \ The stable and the unstable manifold of a
perturbed saddle.}

\medskip

With respect to the two scenarios of Fig. 1, there has been a
(somewhat incorrect) perception that the scenario depicted in Fig.
1(a) is much more interesting than that in Fig. 1(b), mainly due
to the fact that Fig. 1(a) induces complicated dynamics while Fig
1(b) appears simple. A computational procedure (the Melnikov's
method) has been developed to distinguish the two scenarios. See
Sects. 4.5 and 4.6 of \cite{GH} for an exposition on Melnikov's
method and its applications to a list of periodically perturbed
equations. Historically, homoclinic tangles of Fig. 1(a) have been
a strong influence, yet until very recently in the analysis of
concrete equations, little has been rigorously proved beyond the
existence of embedded horseshoes.

In \cite{WO} and \cite{WOk}, we have connected the studies of
periodically perturbed equations in the vicinity of a dissipative
homoclinic solution to the studies of two specific class of maps.
The main idea, motivated by \cite{AS}, is to abandon the globally
defined time-T maps, and in its stead to compute the return maps
in the extended phase space around the unperturbed homoclinic
loop. It has turned out that (1) the dynamics of homoclinic
tangles of Fig. 1(a) are those of a family of infinitely wrapped
horseshoe maps \cite{WOk} (see Sect. \ref{s1.2} for details), and
(2) the dynamics of Fig. 1(b) are those of a family of rank one
maps \cite{WO}. With the return maps obtained in \cite{WO} and in
\cite{WOk}, we are now able to apply many existing theories on
maps to the studies of periodically perturbed equations with
dissipative homoclinic saddle. Among the theories directly applied
are the Newhouse theory \cite{N}, \cite{PT} on homoclinic
tangency; the theory of SRB measures \cite{Si}, \cite{R},
\cite{Bo}; the theory of H\'enon-like attractors \cite{BC2},
\cite{MV}, \cite{BY}; and the recent theory of rank one chaos
\cite{WY1}-\cite{WY3} based on the theory of Benedicks and
Carleson on strongly dissipative H\'enon maps \cite{BC2}. This
method has led us to many new results.

In this paper we illustrate how to apply the analysis of \cite{WO}
and \cite{WOk} to the studies of concrete equations. We use the
equation
\begin{equation}\label{f1-abs}
\frac{d^2 q}{dt^2} + (\lambda - \gamma q^2) \frac{d q}{dt} - q +
q^3 = \mu q^2 \sin \omega t
\end{equation}
as a demonstrative example to prove, among many other things, that
there exist positive measure sets of parameters $(\lambda, \gamma,
\mu, \omega)$ corresponding respectively to the dynamics scenarios
of Fig. 1(a) and Fig. 1(b), so that equation (\ref{f1-abs}) admits
strange attractors with SRB measures. We note that the unperturbed
part of equation (\ref{f1-abs}) has been studied previously in
\cite{HR}.

\section{Statement of results}\label{s1}
Let us re-write equation (\ref{f1-abs}) as
\begin{equation}\label{f1-s1}
\begin{split}
\frac{d q}{dt} & = p \\
\frac{d p}{dt} & = - (\lambda - \gamma q^2)p + q -
q^3 + \mu q^2 \sin \theta \\
\frac{d \theta}{dt} & = \omega.
\end{split}
\end{equation}
The space of $(q, p, \theta)$ is the extended phase space for
equation (\ref{f1-abs}) where $\theta \in S^1$.

\medskip

\subsection{Return maps in extended phase space}\label{s1.1}
We call a saddle fixed point of a two dimensional autonomous
system a {\it dissipative} saddle if the magnitude of the
associated negative eigenvalue is larger than that of the positive
eigenvalue. We call a saddle point a {\it homoclinic} saddle if it
is approached by one solution, which we call a {\it homoclinic}
solution, from both the positive and the negative directions of
time. The autonomous part of equation (\ref{f1-abs}) is written as
\begin{equation}\label{f1-s1a}
\frac{d^2 q}{dt^2} + (\lambda - \gamma q^2) \frac{d q}{dt} - q +
q^3 = 0,
\end{equation}
where $\lambda > 0$ and $\gamma$ are parameters. Let $(q, p) \in
{\mathbb R}^2$ be the phase space where $p = \frac{dq}{dt}$. We
start with a previously existed result on (\ref{f1-s1a})
\cite{HR}.
\begin{prop}[{\bf Dissipative homoclinic saddle}]\label{prop1}
There exists $\lambda_0 > 0$ sufficiently small, such that for
$\lambda \in [0, \lambda_0)$, there exists a $\gamma_{\lambda},
|\gamma_{\lambda}| < 10 \lambda$ such that for $\gamma
=\gamma_{\lambda}$;

(i) equation (\ref{f1-s1a}) has a homoclinic solution for $(q, p)
= (0, 0)$, which we denote as
$$
\ell_{\lambda} = \{ \ell_{\lambda}(t) = (a_{\lambda}(t),
b_{\lambda}(t)), \ t \in {\mathbb R} \},
$$
satisfying $a_{\lambda}(t) > 0$ for all $t$;

 (ii) for any given $L > 0$, there exists a $K(L)$ independent of
$\lambda$, such that for all $t \in [-L, L]$,
$$
|\ell_{\lambda}(t) - \ell_0(t)| < K(L) \lambda.
$$
\end{prop}

See also \cite{W} for an alternative proof of this proposition.
Note that $\ell_0(t)$ is a homoclinic solution of equation
(\ref{f1-s1a}) for $\lambda = \gamma = 0$. Let
\begin{equation}\label{alpha-beta}
\alpha = \frac{1}{2}(\sqrt{\lambda^2 + 4} + \lambda), \ \ \ \beta
= \frac{1}{2}(\sqrt{\lambda^2 + 4} - \lambda).
\end{equation}
Then $- \alpha, \beta$ are the two eigenvalues of  $(0, 0)$. Since
$-\alpha + \beta = -\lambda < 0$, $(0, 0)$ is a dissipative saddle
provided that $\lambda > 0$.

Let $\lambda_0 >0$ be sufficiently small. We regard the value of
$\lambda_0$ as been fixed. We say that $\lambda \in (0,
\lambda_0)$ is {\it non-resonant} if $\alpha, \beta$ satisfy
certain Diophantine non-resonance condition. Namely,

\smallskip

{\bf (H)} \ {\it there exists $d_1, d_2> 0$ so that for all $n, m
\in {\mathbb Z}^+$,
$$
|n \alpha - m \beta| > d_1 (|n| + |m|)^{-d_2}.
$$
}

\smallskip

From this point on we fixed the value of $\lambda \in (0,
\lambda_0)$ and assume that it satisfies (H). Let
$\gamma_{\lambda}$ be as in Proposition \ref{prop1}. We introduce
a new parameter $\rho$ for $\gamma$ by letting
\begin{equation}
\gamma = \gamma_{\lambda} + \mu \rho.
\end{equation}
$\omega, \rho, \mu$ are now the parameters of equation
(\ref{f1-abs}).

\medskip

We turn now to equation (\ref{f1-s1}).  Let $\ell =
\ell_{\lambda}$ be the homoclinic loop of Proposition \ref{prop1}
in the space of $(q, p)$. We construct a small neighborhood of
$\ell$ by taking the union of a small neighborhood
$U_{\varepsilon}$ of $(0, 0)$ and a small neighborhood $D$ around
$\ell$ out of $U_{\frac{1}{4} \varepsilon}$. See Fig. 2. Let
$\sigma^{\pm} \in U_{\varepsilon} \cap D$ be the two line segments
depicted in Fig. 2, both perpendicular to the homoclinic solution.
In the space of $(q, p, \theta)$ we denote
$$
{\mathcal U}_{\varepsilon} = U_{\varepsilon} \times S^1, \ \ \
{\bf D} = D \times S^1
$$
and let
$$
\Sigma^{\pm} = \sigma^{\pm} \times S^1.
$$

\begin{picture}(7, 5.5)
\put(4,0){ \psfig{figure=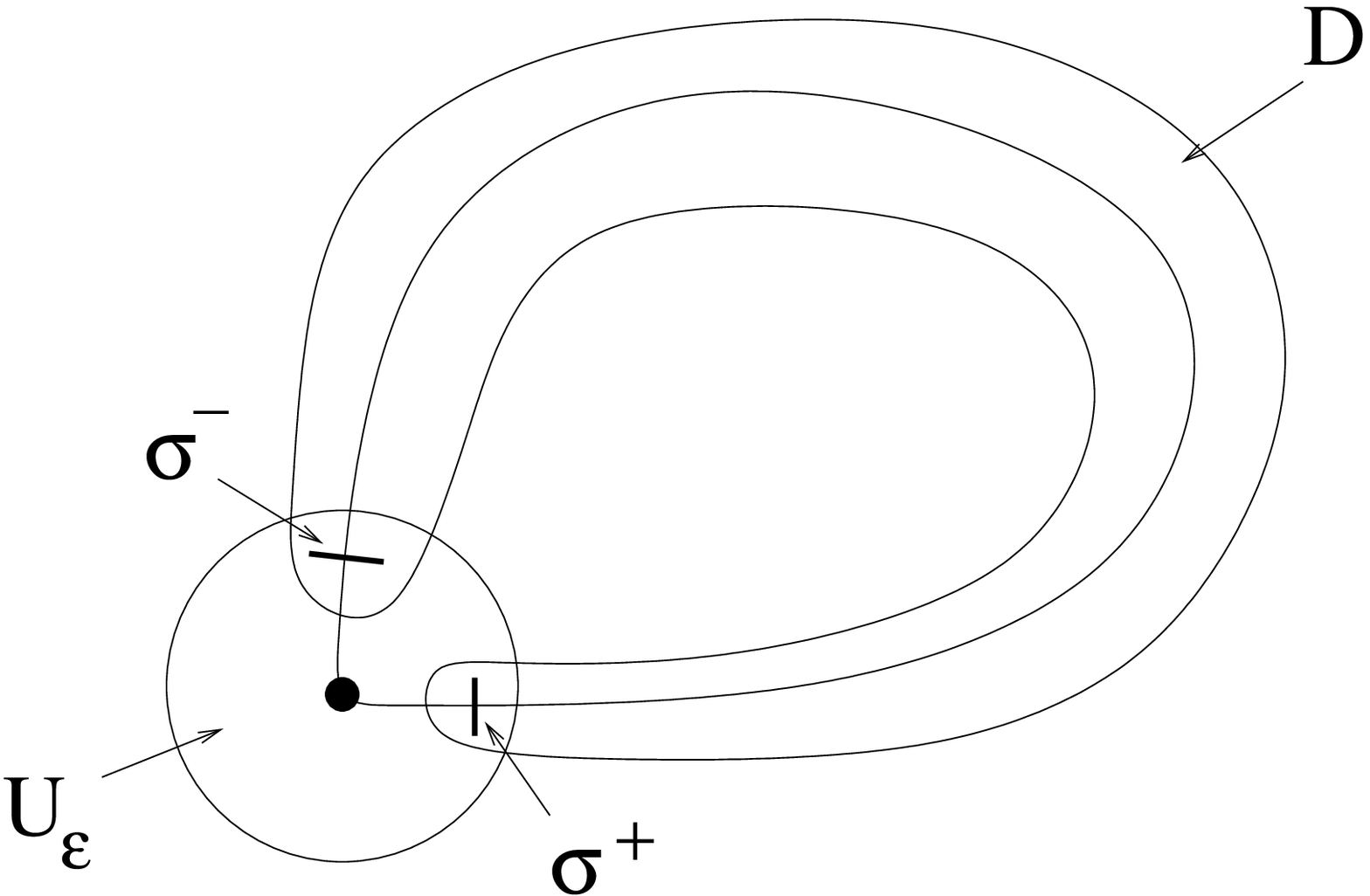,height = 5cm, width =
7cm} }
\end{picture}

\medskip

\centerline{Fig. 2 \ $U_{\varepsilon}$, $D$ and $\sigma^{\pm}$.}

\medskip

Let ${\mathcal N}: \Sigma^+ \to \Sigma^-$ be the maps induced by
the solutions of (\ref{f1-s1}) on ${\mathcal U}_{\varepsilon}$ and
${\mathcal M}: \Sigma^- \to \Sigma^+$ be the maps induced by the
solutions of (\ref{f1-s1}) on ${\bf D}$. See Fig. 3. Let
${\mathcal F}_{\omega, \rho, \mu}:= {\mathcal N} \circ {\mathcal
M}: \Sigma^- \to \Sigma^-$ be the return map defined by equation
(\ref{f1-s1}). We compute ${\mathcal F} = {\mathcal F}_{\omega,
\rho, \mu}$ following the steps of \cite{WO} and \cite{WOk}.

\smallskip

\begin{picture}(8, 6)
\put(3.5,0){ \psfig{figure=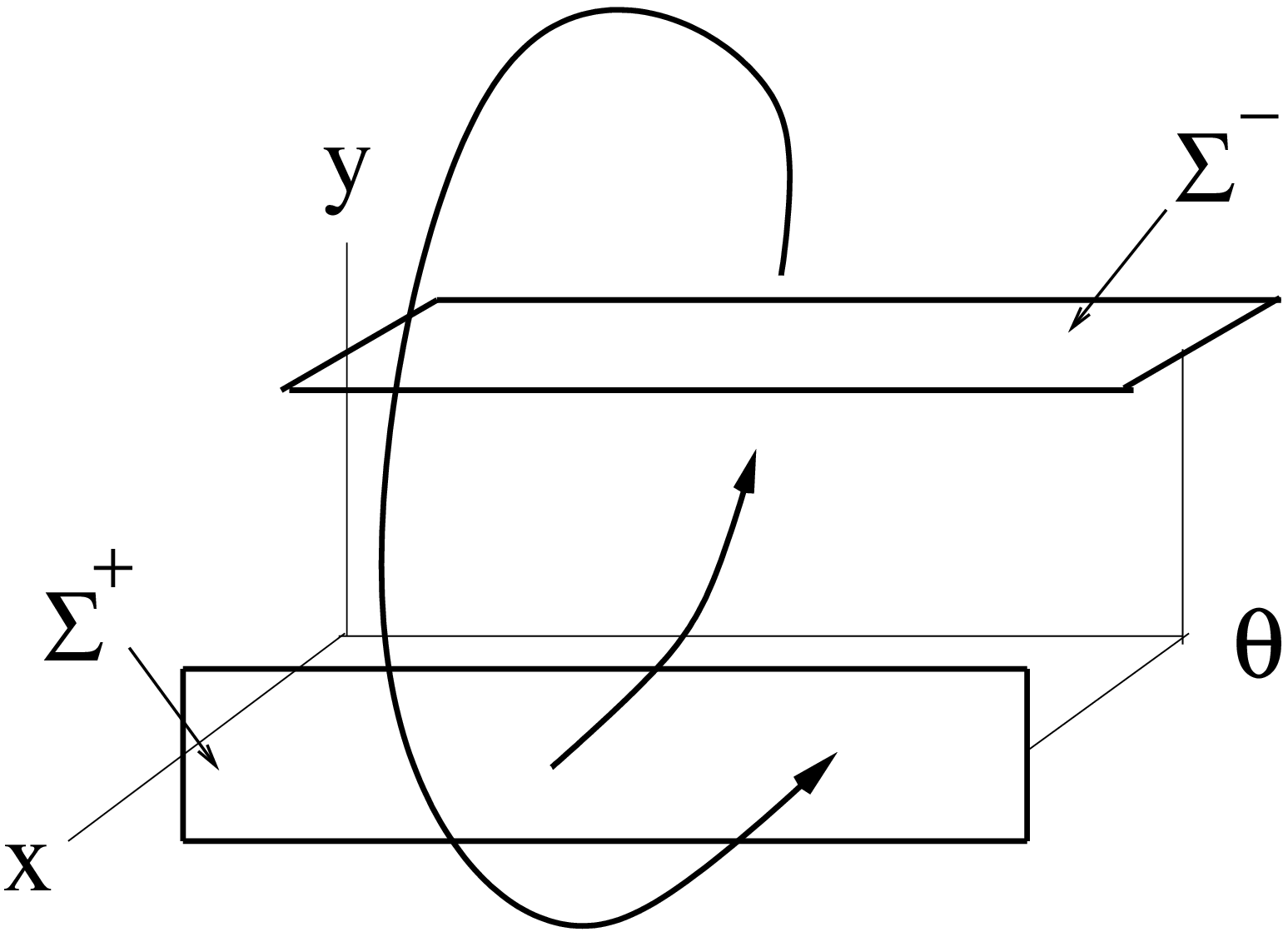,height = 5cm, width =
8cm} }
\end{picture}

\medskip

\centerline{Fig. 3 \ ${\mathcal N}$ and ${\mathcal M}$.}

\medskip

We now define the range of parameters for this paper. Let
$R_{\omega}$ be sufficiently large, and $R_{\omega}, R_{\rho}$ and
$R_{\mu}$ be such that
$$
R_{\omega} << R_{\rho} << R_{\mu}.
$$
We only consider parameters $(\omega, \rho, \mu)$ inside of
$$
{\mathbb P} = \{ (\omega, \rho, \mu): \ \omega \in (0,
R_{\omega}), \ \rho \in (0, R_{\rho}), \ \mu \in (0, R_{\mu}^{-1})
\}.
$$

In the extended phase space $(q, p, \theta)$, $q = p = 0$ is a
periodic orbit. We study equation (\ref{f1-s1}), denoting the
stable and the unstable manifold of the solution $(q, p) = (0, 0)$
as $W^s$ and $W^u$ respectively. First we have

\begin{theorem}[{\bf Surface of tangency}]\label{th1}
There exists a continuous function $S^*(\omega, \mu): (0,
R_{\omega}) \times (0, R_{\mu}^{-1}) \to (0, R_{\rho})$ such that

(a) $W^u \cap W^s \neq \emptyset$ if $0 < \rho < S^*_{\mu}(\omega,
\mu)$;

(b) $W^u \cap W^s = \emptyset$ if $S^*_{\mu}(\omega, \mu) < \rho <
R_{\rho}$.
\end{theorem}

So under the surface $S^*$ defined by $\rho  = S^*(\omega, \mu)$
is the scenario of homoclinic tangles of Fig. 1(a) and above this
surface is the scenario of Fig. 1(b). $S^*$ is one of the surfaces
(the lowest one of the three) depicted in Fig. 6 in Sect.
\ref{s1.3}.

\subsection{Homoclinic tangles}\label{s1.2}
Assume that $(\omega, \rho, \mu) \in {\mathbb P}$ is beneath the
surface $S^*$ of Theorem \ref{th1}. We focus on the set of
solutions that stay inside of ${\mathcal U}_{\varepsilon} \cup
{\bf D}$ for all time. This set of solutions is the homoclinic
tangle in the vicinity of $\ell_{\lambda}$. If $(\omega, \rho,
\mu) \in {\mathbb P}$ is beneath $S^*$, then the return map
${\mathcal F} = {\mathcal F}_{\omega, \rho, \mu}$ is only
partially defined on $\Sigma^-$ because some of the images of
$\Sigma^-$ under ${\mathcal M}$ would first hit on the wrong side
of the local stable manifold of $(q, p) = (0, 0)$ in $\Sigma^+$,
then sneak out of ${\mathcal U}_{\varepsilon}$. See Fig. 4.

\begin{picture}(9, 5.5)
\put(4.5,0){ \psfig{figure=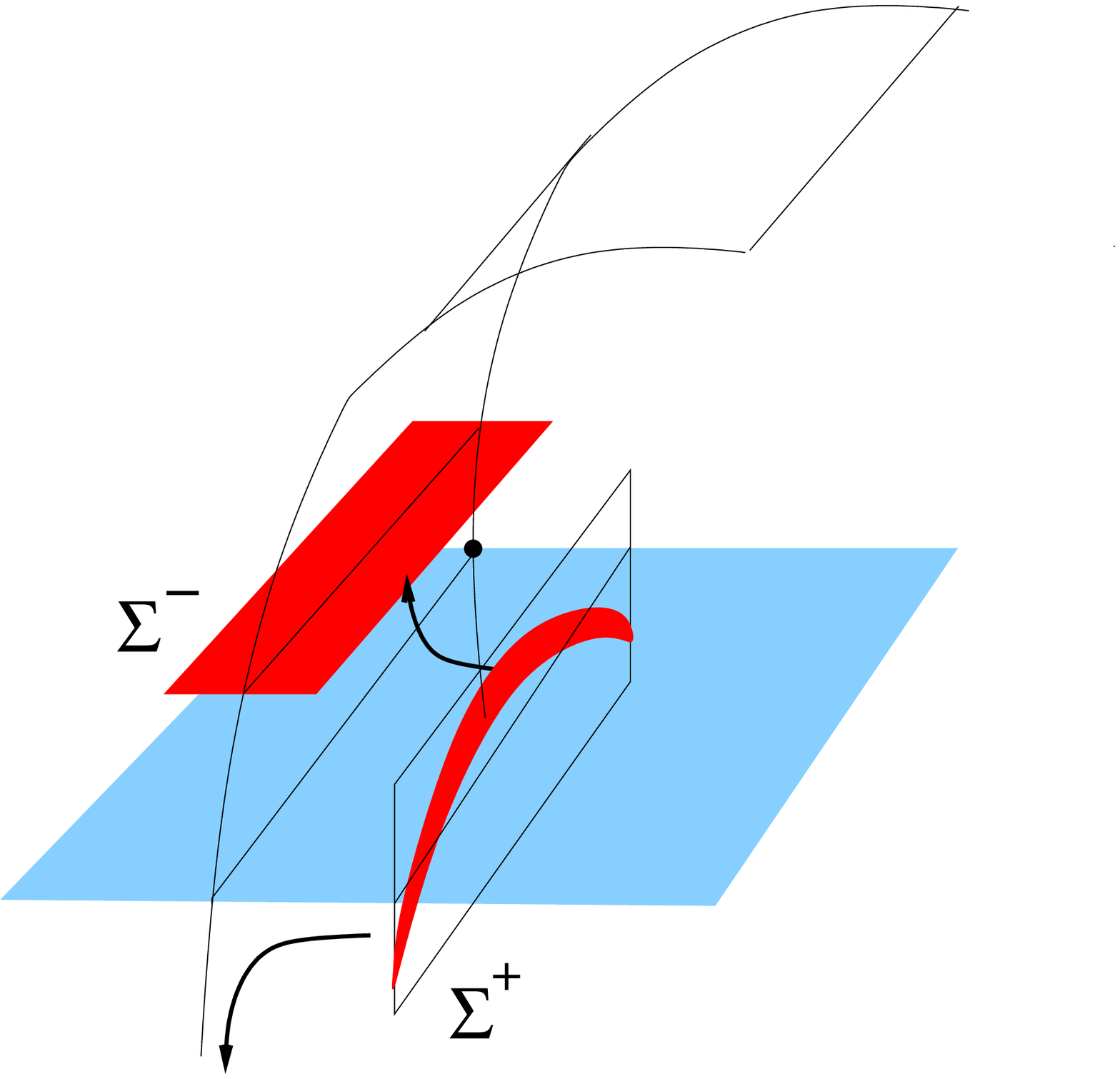,height = 5cm, width =
8cm} }
\end{picture}

\medskip

\centerline{Fig. 4 \  Partial returns to $\Sigma^-$.}

\medskip

We prove that ${\mathcal F}$ is an infinitely wrapped horseshoe
map, the geometric structure of which is as follows. Take an
annulus ${\mathcal A} = S^1 \times I$ (This is $\Sigma^-$ for
${\mathcal F}$). We represent points in $S^1$ and $I$ by using
variables $\theta$ and $z$ respectively. We call the direction of
$\theta$ the horizontal direction and the direction of $z$ the
vertical direction. To form an infinitely wrapped horseshoe map,
which we denote as ${\mathcal F}$, we first divide ${\mathcal A}$
into two vertical strips, denoted as $V$ and $U$. ${\mathcal F}: V
\to {\mathcal A}$ is defined on $V$ but not on $U$. We compress
$V$ in the vertical direction and stretch it in the horizontal
direction, making the image infinitely long towards both ends.
Then we fold it and wrap it around the annulus ${\mathcal A}$
infinitely many times. See Fig. 5.

\begin{picture}(7, 5.5)
\put(3.7,0){ \psfig{figure=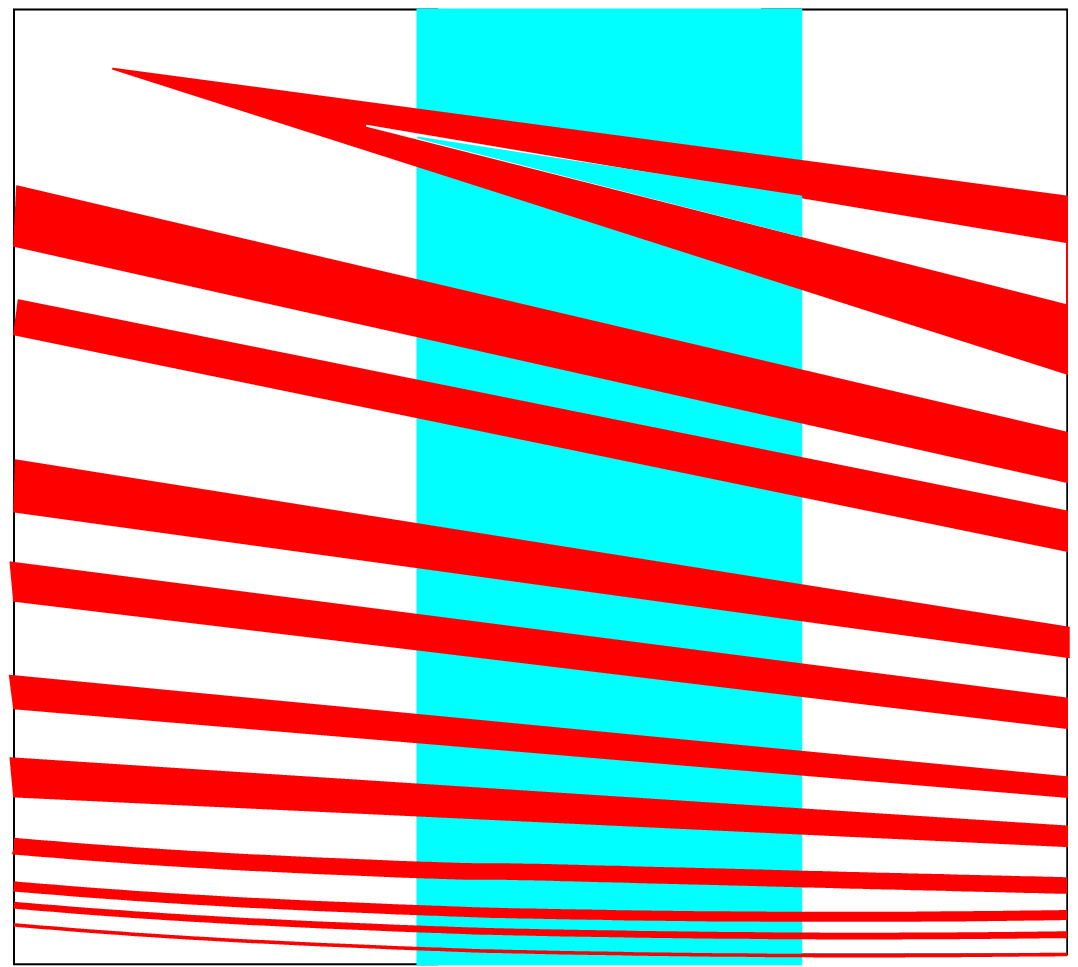,height = 5cm, width =
7cm} }
\end{picture}

\medskip

\centerline{Fig. 5 \ Infinitely wrapped horseshoe maps.}

\medskip

Let
$$
\Omega_{\mathcal F} = \{ (\theta, z) \in V: \ {\mathcal
F}^n(\theta, z) \in V \ \ \ \forall n \geq 0 \}, \ \ \ \ \ \ \
\Lambda_{\mathcal F} = \cap_{n\geq 0} {\mathcal
F}^n(\Omega_{\mathcal F}).
$$
The dynamic structure of the homolinic tangle in the vicinity of
$\ell_{\lambda}$ is manifested in that of $\Lambda_{\mathcal F}$.
$\Lambda_{\mathcal F}$ obviously contain a horseshoe of infinitely
many symbols. This horseshoe covers Smale's horseshoe and all its
variations. It is the one that resides inside all homoclinic
tangles.

The structure of $\Omega_{\mathcal F}$ and $\Lambda_{\mathcal F}$
depend sensitively on the location of the folded part of
${\mathcal F}(V)$. If this part is deep inside of $U$, then the
entire homoclinic tangle is reduced to one horseshoe of infinitely
many symbols. If it is located inside of $V$, then the homoclinic
tangles are likely to have attracting periodic solutions, or sinks
and observable chaos associated with non-degenerate transversal
homoclinic tangency. In particular, we have

\begin{theorem}[{\bf Dynamics of homoclinic tangles}]\label{th2}
There exists an open domain in the plane of $(\omega, \rho)$,
which we denote as $DH \subset (0, R_{\omega}) \times (0,
R_{\rho})$, such that $DH \times (0, R_{\mu}^{-1})$ is strictly
under the surface $S^*$ in ${\mathbb P}$; and for any fixed
$(\omega, \rho) \in DH$, the return maps ${\mathcal F}_{\mu} =
{\mathcal F}_{\omega, \rho, \mu}$, as a one parameter family in
$\mu$, satisfy the follows:

(i) There exists a sequence of $\mu$, accumulating at $\mu = 0$,
which we denote as
$$
R_{\mu}^{-1} \geq \mu_1^{(r)} > \mu_1^{(l)} > \cdots > \mu_n^{(r)}
> \mu_n^{(l)} > \cdots > 0
$$
such that for all $\mu \in [\mu^{(l)}_n, \mu^{(r)}_n]$, ${\mathcal
F} = {\mathcal F}_{\mu}: \Lambda_{\mathcal F} \to
\Lambda_{\mathcal F}$ conjugates to a full shift of countably many
symbols.

(ii) Inside each of the intervals $[\mu^{(r)}_{n+1},
\mu^{(l)}_n]$, $n > 0$,  there exists  $\mu$ so that ${\mathcal
F}_{\mu}$ admits stable periodic solutions.

(iii) There are also parameters in each of the intervals
$[\mu^{(r)}_{n+1}, \mu^{(l)}_n]$, $n > 0$, such that ${\mathcal
F}_{\mu}$ admits non-degenerate transversal homoclinic tangency of
a dissipative saddle fixed point.
\end{theorem}

\medskip

\noindent {\bf Remarks:} \ (1) All results stated so far in this
subsection about homoclinic tangles are new. We are not aware of
any previous results of the same kind in the analysis of any
concrete time-periodic second order equations.

(2) Theorem \ref{th2}(i) is much stronger than the mere existence
of an embedded horseshoe. For these parameters, the entire
homoclinic tangle around $\ell_{\lambda}$ contains not only
nothing less \cite{S}, but also nothing more, than a uniformly
hyperbolic invariant set conjugating to a horseshoe of countably
many symbols. Since the attractive basin of such horseshoe is of
Lebesgue measure zero, these tangles are not observable if one
adopt the probabilistic point of view that only events of positive
Lebesgue measures are observable in phase space.

(3) Sinks of Theorem \ref{th2}(ii) are sinks of very strong
contraction. They are not related to those obtained through
Newhouse theory. The existence of Newhouse sinks follows from
Theorem \ref{th2}(iii).

(4) It also follows from combining Theorem \ref{th2}(iii) with the
result of \cite{MV} that there are positive measure set of
parameters $\mu$, such that the return maps ${\mathcal F}$ admit
H\'enon-like attractors. In contrast to the horseshoes of Theorem
\ref{th2}(i), these are chaos in observable forms (SRB measures)
\cite{MV}, \cite{BY}.

\subsection{Invariant curves and rank one chaos}\label{s1.3}
Assume that $(\omega, \rho, \mu) \in {\mathbb P}$ is above the
surface $S^*$ of Theorem \ref{th1}. We are in the case of Fig.
1(b). The return maps ${\mathcal F}$ are now well-defined on
$\Sigma^-$, and they are not in anyway less interesting than those
for homoclinic tangles of Section \ref{s1.2} in terms of chaos
dynamics. For these parameters, not only the return maps are rank
one maps, to which the theory of \cite{WY1} and \cite{WY2}
directly apply, but also they fall naturally into the specific
category of rank one maps previously studied by Wang and Young in
\cite{WY4}. For these parameters the dynamics of equation
(\ref{f1-s1}) around $\ell_{\lambda}$ are determined by the
magnitude of the forcing frequency $\omega$. When the forcing
frequency $\omega$ is small, we obtain, for all $\mu$ sufficiently
small, an attracting tori in the extended phase space. In
particular, we obtain an attracting tori consisting of
quasi-periodic solutions for a set of $\mu$ with positive Lebesgue
density at $\mu = 0$. As $\omega$ increases, the attracting tori
is dis-integrated into isolated periodic sinks and saddles.
Increasing the forcing frequency $\omega$ further, the stable and
the unstable manifold of these periodic saddles will fold, and
intersect to create horseshoe and strange attractors. We have, in
particular, that these are strange attractors with SRB measures.
In particular, we have
\begin{theorem}[{\bf Invariant curves and quasi-periodic tori}]\label{th3}
There exists a $\rho_0>0$ and a continuous function $Q: (\rho_0,
R_{\rho}) \times (0, R_{\mu}^{-1}) \to (0, R_{\omega})$, such that
for any given $(\omega, \rho, \mu) \in {\mathbb P}$ that is on the
left hand side of the surface $Q$ defined by $\omega = Q(\rho,
\mu)$, that is, these satisfying $\omega < Q(\rho, \mu)$, we have
the follows:

(1) The return map ${\mathcal F}: \Sigma^- \to \Sigma^-$ admits a
globally attracting, simply closed invariant curve; and the maps
induced on this invariant curve is equivalent to a circle
diffeomorphism.

(2) There is an open set $DI \subset (0, R_{\omega}) \times
(\rho_0, R_{\rho})$ in the $(\omega, \rho)$-plane, such that $DI
\times (0, R_{\mu}^{-1})$ is located on the left hand side of the
surface $Q$. Let $(\omega, \rho) \in DI$ be fixed and regard
${\mathcal F}_{\mu} = {\mathcal F}_{\omega, \rho, \mu}$ as a one
parameter family. Then there exists a set of $\mu$ of positive
Lebesgue density at $\mu = 0$, such that the circle diffeomorphism
induced on the attracting invariant curve above is with an
irrational rotation number.
\end{theorem}

\noindent {\bf Remark:} \ See Fig. 6 for $Q$ inside of ${\mathbb
P}$. Similar results are previously obtained for periodically
kicked limit cycles in \cite{WY4}, \cite{WY5}, \cite{LWY},  and
for periodically perturbed hyperbolic periodic solutions in
\cite{Lev2}, \cite{H}.

\smallskip

We turn to the case of strange attractors and rank one chaos.
\begin{theorem}[{\bf Strange attractors and rank one chaos}]\label{th4}
Let $S^*(\omega, \mu)$ be as in Theorem \ref{th2} and
$$
S(\omega, \mu) = (1 + \omega^\frac{1}{2}) S^*(\omega, \mu).
$$
Let $(\omega, \rho, \mu) \in {\mathbb P}$ be such that
$$
S^*(\omega, \mu) < \rho < S(\omega, \mu).
$$
Then,

(1) the return map ${\mathcal F}: \Sigma^- \to \Sigma^-$ admits a
global attractor that is strange in the sense that it contains an
embedded horseshoe provided that $\omega > 100$; and

(2) there exists an open domain $DC \subset (100, R_{\omega})
\times (0, R_{\rho})$ in the $(\omega, \rho)$-plane, such that $DC
\times (0, R_{\mu}^{-1})$ is located in between $S^*$ and the
surface $S$ defined by $\rho = S(\omega, \mu)$. Let $(\omega,
\rho) \in DC$ be fixed and regard ${\mathcal F}_{\mu} = {\mathcal
F}_{\omega, \rho, \mu}$ as a one parameter family. Then there
exists a set of $\mu$ of positive Lebesgue density at $\mu = 0$,
such that ${\mathcal F}_{\mu}$ admits a strange attractor with an
ergodic SRB measure $\nu$. Furthermore, almost every point of
$\Sigma^-$ is generic with respect to $\nu$.
\end{theorem}

\noindent {\bf Remarks:} \ (1) Various dynamical scenarios of
Theorems \ref{th1}-\ref{th4} are put together in Fig. 6. This
figure is self-explanatory. Results similar to Theorem  \ref{th4}
are previously obtained for periodically kicked limit cycles in
\cite{WY4}, \cite{WY5}, \cite{LWY}, and for certain slow-fast
systems in \cite{GWY}.

\begin{picture}(9, 7.5)
\put(2,0){ \psfig{figure=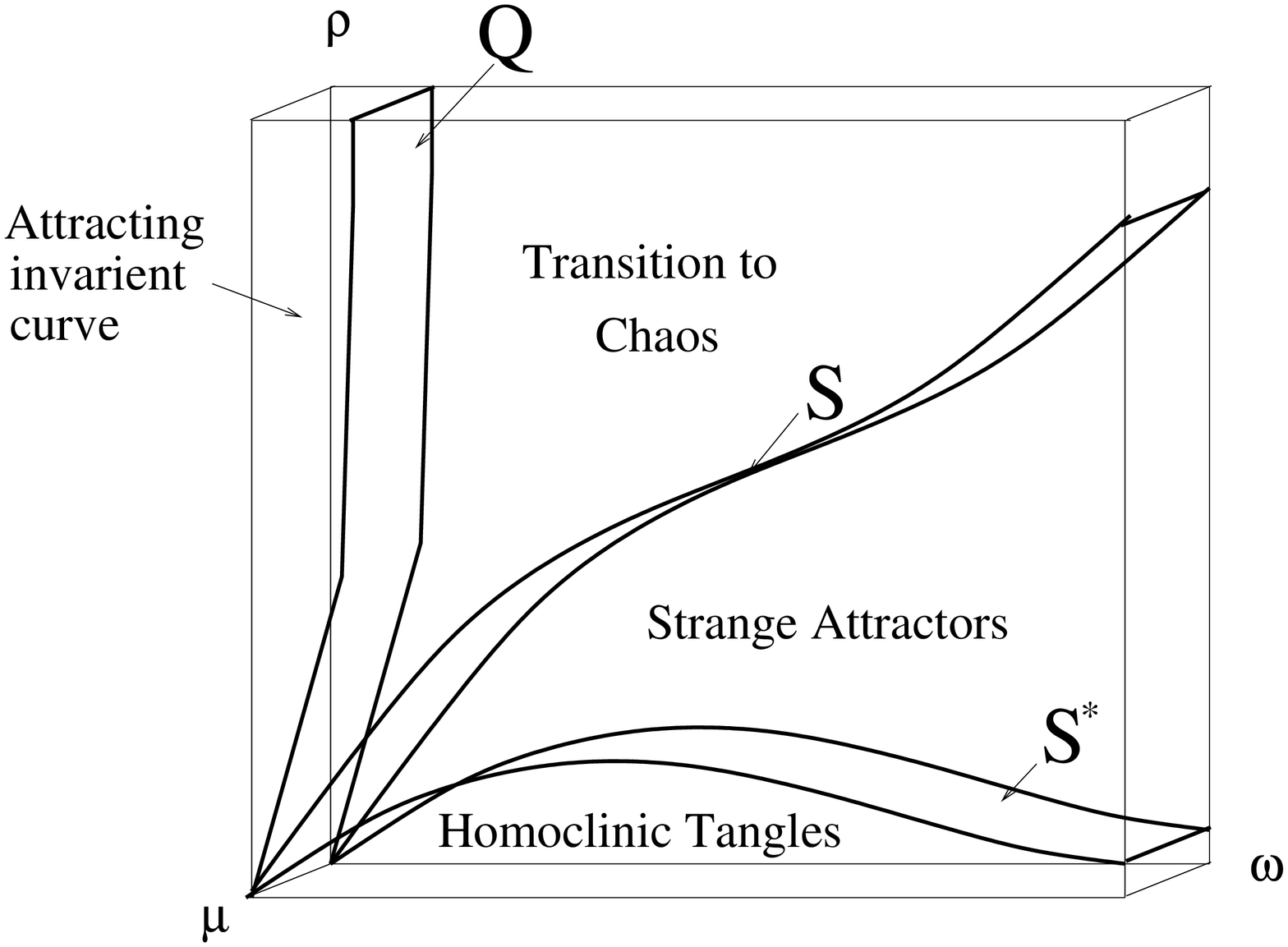,height = 7.5cm, width =
8cm} }
\end{picture}

\vskip .1in

\centerline{Fig. 6 \ Dynamical scenarios in parameter space.}

\medskip

(2) Theorems \ref{th3} and \ref{th4} remain valid if we change the
forcing function $\mu q^2 \sin \omega t$ in (\ref{f1-abs}) to
simply $\mu \sin \omega t$. This case is covered by \cite{WO}.

(3) SRB measures are constructed initially for uniformly
hyperbolic systems by Sinai, Ruelle and Bowen (\cite{Si},
\cite{R}, \cite{Bo}). They are observable objects representing
statistical order in chaos. We refer to \cite{Y} for a review on
the theory of SRB measures. \cite{WY1} contains not only the
existence of SRB measures but also a comprehensive dynamical
profile for the good maps, including geometric structures of the
attractors and statistical properties. We have opted to limit our
statement to SRB measures, but all aspects of that larger
dynamical picture in fact apply.

(4) Finally we note that the open domains $DH, DI$ and $DC$ in
Theorems \ref{th2}-\ref{th4} are not in any sense small in size.
One main restriction is that their respective product to $(0,
R_{\mu}^{-1})$ do not intersect $Q$, $S$ and $S^*$. This is
essentially the only restriction for $DI$ of Theorem \ref{th3}.
For $DH$ and $DC$ we also require that $\omega$ is reasonably
large.

\section{Derivation of the return maps}\label{s2}
In this section we derive the return maps ${\mathcal F} =
{\mathcal N} \circ {\mathcal M}$ of Sect. \ref{s1.1} for equation
(\ref{f1-s1}) assuming Proposition \ref{prop1}. In Sect.
\ref{s2.1} we introduce coordinate changes to transform the
equations into canonical forms, which we will use in Sect.
\ref{s2.2} to compute the return maps.

\smallskip

\subsection{Equation in canonical forms}\label{s2.1} In Sect.
\ref{s2.1}A we normalize the linear part. In Sect. \ref{s2.1}B we
linearize equation (\ref{f1-s1}) on ${\mathcal U}_{\varepsilon}$,
the $\varepsilon$-neighborhood of the saddle point. In Sect.
\ref{s2.1}C we introduce a set of new variables and derive the
corresponding equations on ${\bf D}$, a small neighborhood around
the homoclinic solution $\ell_{\lambda}$ out of ${\mathcal
U}_{\frac{1}{4} \varepsilon}$.

\medskip

\noindent {\bf A. Normalizing the linear part} \ Let $\lambda \in
(0, \lambda_0)$ be fixed satisfying (H) in Sect. \ref{s1.1}, and
$\gamma = \gamma_{\lambda} + \mu \rho$ where $\gamma_{\lambda}$ is
as in Proposition \ref{prop1}. We re-write equation (\ref{f1-s1})
as

\begin{equation}\label{f1-s2.1a}
\begin{split}
\frac{d q}{dt} & = p \\
\frac{d p}{dt} & = - (\lambda - \gamma_{\lambda} q^2)p + q -
q^3 + \mu (\rho q^2 p +  q^2 \sin \theta) \\
\frac{d \theta}{dt} & = \omega.
\end{split}
\end{equation}

To put the linear part of equation (\ref{f1-s2.1a}) in canonical
form, we introduce new variables $(x, y)$ so that
\begin{equation}\label{f2-s2.1a}
q = x + \alpha y, \ \ \ \ p = -\alpha x + y,
\end{equation}
where $\alpha = \frac{1}{2}(\lambda + \sqrt{\lambda^2 + 4})$ is as
in (\ref{alpha-beta}) in Sect. \ref{s1.1}. In reverse we have
\begin{equation}\label{f3-s2.1a}
x = \frac{1}{1 + \alpha^2}( q - \alpha p) \ \ \ \ y = \frac{1}{1+
\alpha^2} (\alpha q + p).
\end{equation}
The new equations for $(x, y)$ are
\begin{equation}\label{f4-s2.1a}
\begin{split}
\frac{d x}{dt} & = - \alpha x + f(x, y) + \mu (A(x, y) \rho + C(x, y) \sin \theta) \\
\frac{d y}{dt} & = \beta y + g(x, y) + \mu(B(x, y) \rho + D(x, y) \sin \theta )   \\
\frac{d \theta}{dt} & = \omega
\end{split}
\end{equation}
where $\beta = \alpha^{-1}$ is again as in (\ref{alpha-beta}), and
\begin{equation*}
\begin{split}
f(x, y) & = \frac{\alpha}{1 + \alpha^2}\left(\gamma_{\lambda}(x +
\alpha y)^2(y -
\alpha x) + (x + \alpha y)^3\right),\\
g(x, y) & = \frac{-1}{1 + \alpha^2}\left(\gamma_{\lambda}(x +
\alpha y)^2(y - \alpha x) + (x + \alpha y)^3\right);
\end{split}
\end{equation*}
\begin{equation*}
\begin{split}
A(x, y) & = \frac{\alpha}{1 + \alpha^2}(x + \alpha y)^2(y -
\alpha x), \ \ \ \ B(x, y) = \frac{-1}{1 + \alpha^2} (x + \alpha y)^2(y - \alpha x); \\
C(x, y) & = \frac{\alpha}{1 + \alpha^2}(x + \alpha y)^2, \ \ \ \
D(x, y) = \frac{-1}{1 + \alpha^2}(x + \alpha y)^2.
\end{split}
\end{equation*}

Observe that the functions $f, g, A, B, C, D$ are all functions of
$x, y$ of order at least two. In comparison to the equations
studied in \cite{WOk}, equation (\ref{f4-s2.1a}) is in a slightly
different form. In the rest of this section we present the
derivations of the return maps following mainly \cite{WOk}. To
avoid repetitive writings, we will refer the reader to the
corresponding part of \cite{WOk} for the details of a few
technical proofs along the way.

\smallskip

\noindent {\bf Notations} \ $(\omega, \rho, \mu)$ are the three
parameters of equation (\ref{f1-s1}). In what follows we replace
$\mu$ by $\ln \mu$. $\mu \in (0, R_{\mu}^{-1})$ corresponds to
$\ln \mu \in (-\infty, - \ln R_{\mu})$. Denote $p = (\omega, \rho,
\ln \mu)$, $p \in {\mathbb P} =  (0, R_{\omega}) \times (0,
R_{\rho}) \times (-\infty, - \ln R_{\mu})$. $\varepsilon$ is not a
parameter of equation (\ref{f1-s1}) but an auxiliary parameter
representing the size of $U_{\varepsilon}$. In what follows
$\varepsilon$ is a positive number sufficiently small satisfying
\begin{equation}\label{f5-s2.1a}
R_{\mu} >> \varepsilon^{-1} >> R_{\rho} >> R_{\omega}.
\end{equation}

The letter $K$ is reserved as a generic constant, the value of
which varies from line to line. $K$ is allowed to depend on
parameters $\omega, \rho$ and $\varepsilon$ but not $\mu$. We also
use $K(\varepsilon)$ to make explicit when a constant is a
dependent of $\varepsilon$. Letter $K$ standing alone represents a
constant independent of both $\varepsilon$ and $\mu$.

The intended formula for the return maps would contain explicit
terms and ``error" terms, and we aim on $C^3$-control on all error
terms in this paper. To facilitate our presentation we adopt a
specific conventions for indicating controls on magnitude. For a
given constant, we write ${\mathcal O}(1)$, ${\mathcal
O}(\varepsilon)$ or ${\mathcal O}(\mu)$ to indicate that the
magnitude of the constant is bounded by $K$, $K \varepsilon$ or
$K(\varepsilon) \mu$, respectively. For a function of a set $V$ of
variables and parameters on a specific domain, we write ${\mathcal
O}_V(1), {\mathcal O}_V(\varepsilon)$ or ${\mathcal O}_V(\mu)$ to
indicate that the $C^r$-norm of the function on the specified
domain is bounded by $K, K \varepsilon$ or $K(\varepsilon) \mu$,
respectively. We chose to specify the domain in the surrounding
text rather than explicitly involving it in the notation. For
example, ${\mathcal O}_{\theta, p}(\mu)$ represents a function of
$\theta$, the $C^3$-norm\footnote{We emphasize that $\ln \mu$, not
$\mu$, is regarded as one of the forcing parameters. Derivatives
with respect to $\mu$ and to $\ln \mu$ differ by a factor $\mu$.}
of which with respect to $\theta, \omega, \rho$ and $\ln \mu$ is
bounded above by $K(\varepsilon) \mu$.

\medskip

\noindent {\bf B. Linearization on ${\mathcal U}_{\varepsilon}$.}
\ Let $X, Y$ be such that
\begin{equation}\label{f1-s2.1b}
\begin{split}
x & = X + P(X, Y) + \mu \tilde P(X, Y, \theta; p) \\
y & = Y + Q(X, Y) + \mu \tilde Q(X, Y, \theta; p)
\end{split}
\end{equation}
where $P, Q, \tilde P, \tilde Q$ as functions of $X$ and $Y$ are
real-analytic on $|(X, Y)| < 2 \varepsilon$, and the values of
these functions and their first derivatives with respect to $X$
and $Y$ at $(X, Y) = (0, 0)$ are all zero. As is explicitly
indicated in (\ref{f1-s2.1b}), $P$ and $Q$ are independent of
$\theta$ and $p$. We also assume that
$$
\tilde P(X, Y, \theta+ 2 \pi; p)= \tilde P(X, Y, \theta; p), \ \
\tilde Q(X, Y, \theta+ 2 \pi; p)= \tilde Q(X, Y, \theta; p)
$$
are periodic of period $2 \pi$ in $\theta$ and they are also
real-analytic with respect to $\theta$ and $p$ for all $\theta \in
{\mathbb R}$ and $p \in {\mathbb P}$. We have
\begin{prop}\label{prop1-s2.1b}
There exists a small neighborhood $U_{\varepsilon}$ of $(0, 0)$ in
$(X, Y)$-space, the size of which are completely determined by
equation (\ref{f1-s1}) and $d_1, d_2$ in (H), such that there
exists an analytic coordinate transformation in the form of
(\ref{f1-s2.1b}) on ${\mathcal U}_{\varepsilon} = U_{\varepsilon}
\times S^1$ that transforms equation (\ref{f4-s2.1a}) into
\begin{equation*}
\frac{d X}{dt} = - \alpha X, \ \ \frac{d Y}{dt} = \beta Y, \ \
\frac{d \theta}{dt} = \omega.
\end{equation*}
Moreover, the $C^4$-norms of $P, Q, \tilde P, \tilde Q$ as
functions of $X, Y, \theta, \mu$ are all uniformly bounded from
above by a constant $K$ that is independent of both $\varepsilon$
and $\mu$ on $(X, Y) \in U_{\varepsilon}$, $\theta \in {\mathbb
R}$ and $p \in {\mathbb P}$.
\end{prop}
\noindent {\bf Proof:} \ This is a standard linearization result.
For this proposition to hold we need the non-resonant assumption
(H) and the fact that $f, g, A, B, C, D$ are terms of order higher
than one at $(x, y) = (0, 0)$. See for instance \cite{CLS} for a
proof. \hfill $\square$

\medskip

\noindent {\bf C. A canonical form around homoclinic loop. } \ We
derive a standard form for equation (\ref{f4-s2.1a}) around the
homoclinic loop of equation (\ref{f1-s1a}) outside of ${\mathcal
U}_{\frac{1}{4}\varepsilon}$. Let $\ell = \ell_{\lambda}$ be the
homoclinic solution of Proposition \ref{prop1}. In $(x, y)$-space
we write $\ell$ as
$$
x = a(t), \ \ \ \ \ \ y = b(t)
$$
where $t \in {\mathbb R}$ is the time. Let
$$
(u(t), v(t)) = \left|\frac{d}{dt} \ell(t)\right|^{-1} \frac{d}{dt}
\ell(t)
$$
be the unit tangent vector of $\ell$ at $\ell(t)$. We replace $t$
by $s$ to write this homoclinic loop as $\ell(s) = (a(s), b(s))$.
Let
\begin{equation*}
{\bf e}(s) = (v(s), -u(s)).
\end{equation*}
We now introduce new variables $(s, z)$ through
\begin{equation*}
(x, y) = \ell(s) + z {\bf e}(s).
\end{equation*}
This is to say that
\begin{equation}\label{f1-s2.1c}
x = x(s, z) := a(s) + v(s) z, \ \  y =y(s, z):= b(s) - u(s) z.
\end{equation}
Let us further re-scale $z$ by letting $Z = \mu^{-1} z$. $(s, Z,
\theta)$ are the variables we use in this paragraph.

Before deriving the equations for $(s, Z,\theta)$, we first define
the domain on which these equations are derived. Let $-L^-, L^+$
be the times $\ell = \ell_{\lambda}$ hits $U_{\frac{1}{2}
\varepsilon}$. $L^{\pm}$ are completely determined by
$\ell_{\lambda}$ and $\varepsilon$. Defined
$$
{\bf D} = \{ (s, Z; p): \ s \in [-2L^-, 2L^+], \ |Z| \leq
K_1(\varepsilon), \ p \in {\mathbb P} \}.
$$
for some $K_1(\varepsilon)$ independent of $\mu$. We have
\begin{prop}\label{prop1-s2.1c}
The equation for $(s, Z, \theta)$ on ${\bf D}$ is in the form of
\begin{equation}\label{f2-s2.1c}
\begin{split}
\frac{dZ}{dt} & = E(s) Z + ({\mathbb H}_1(s) \rho+ {\mathbb H_2}(s) \sin \theta)
+ {\mathcal O}_{s, Z, \theta; p}(\mu) \\
\frac{ds}{dt} & = 1 + {\mathcal O}_{s, Z, \theta; p}(\mu) \\
\frac{d \theta}{dt} & = \omega
\end{split}
\end{equation}
where
\begin{equation}\label{f3-s2.1c}
\begin{split}
E(s) & = v^2(s) (-\alpha + \partial_x f(a(s), b(s))) + u^2(s)
(\beta + \partial_y g(a(s), b(s))) \\
& \ \ - u(s) v(s) (\partial_y f(a(s), b(s)) + \partial_x g(a(s),
b(s))),
\end{split}
\end{equation}
and
\begin{equation}\label{f4-s2.1c}
\begin{split}
{\mathbb H}_1(s) & = v(s) A(a(s), b(s)) - u(s) B(a(s), b(s)); \\
{\mathbb H}_2(s) & = v(s) C(a(s), b(s)) - u(s) D(a(s), b(s)).
\end{split}
\end{equation}
\end{prop}
\noindent {\bf Proof:} \ The derivations of equation
(\ref{f2-s2.1c}) is identical to that of equation (2.13) in
\cite{WOk}. We note that $K_1(\varepsilon)$ for ${\bf D}$ is also
defined explicitly in \cite{WOk}. \hfill $\square$

\subsection{Derivation of Return Maps} \label{s2.2}
In this subsection we derive the return map ${\mathcal F} =
{\mathcal N} \circ {\mathcal M}$ by using Propositions
\ref{prop1-s2.1b} and \ref{prop1-s2.1c}. In Sect. \ref{s2.2}A we
defined precisely the surfaces $\Sigma^{\pm}$, and discuss the
issue of coordinate conversions on these surfaces. In Sect.
\ref{s2.2}B we compute ${\mathcal N}$ by using Proposition
\ref{prop1-s2.1b} and in Sect. \ref{s2.2}C we compute ${\mathcal
M}$ by using Proposition \ref{prop1-s2.1c}. ${\mathcal F} =
{\mathcal N} \circ {\mathcal M}$ is then obtained in Sect.
\ref{s2.2}D.

\medskip

\noindent {\bf A. \ Poincar\'e sections $\Sigma^{\pm}$.} \ Let
$(X, Y, \theta)$ be the phase variables of Proposition
\ref{prop1-s2.1b}. We define $\Sigma^{\pm}$ inside of ${\mathcal
U}_{\varepsilon} \cap {\bf D}$ by letting
$$
\Sigma^- = \{ (X, Y, \theta): \ Y = \varepsilon, \ |X| < \mu, \
\theta \in S^1 \},
$$
and
$$
\Sigma^+ = \{ (X, Y, \theta): \ X = \varepsilon, \ |Y| <
K_1(\varepsilon) \mu, \ \theta \in S^1 \}.
$$
$K_1(\varepsilon)$ is the same as in Proposition
\ref{prop1-s2.1c}.

\smallskip

Let $q \in \Sigma^+$ or $\Sigma^-$. We represent $q$ by using the
$(X, Y, \theta)$-coordinates, for which we have $X = \varepsilon$
on $\Sigma^+$ and $Y = \varepsilon$ on $\Sigma^-$. We can also use
$(s, Z, \theta)$ of Proposition \ref{prop1-s2.1c} to represent the
same $q$. Before computing the return maps, we need to first
attend two issues that are technical in nature. First, we need to
derive the defining equations of $\Sigma^{\pm}$ for $(s, Z,
\theta)$. Second, we need to be able to change variables back and
forth from $(X, Y, \theta)$ to $(s, Z, \theta)$ on $\Sigma^{\pm}$.
Let
\begin{equation}
{\mathbb X} = \mu^{-1} X, \ \ \ {\mathbb Y} = \mu^{-1} Y.
\end{equation}
\begin{prop}\label{prop1-s2.2a} (a) On $\Sigma^+$, we have

(i) $s = -L^- + {\mathcal O}_{Z, \theta, p}(\mu)$; and

(ii) $ {\mathbb Y} = (1 + {\mathcal O}(\varepsilon))Z + {\mathcal
O}_{\theta, p}(1) + {\mathcal O}_{Z, \theta, p}(\mu)$.

\noindent (b) On $\Sigma^-$, we have

(i) $s = -L^- + {\mathcal O}_{Z, \theta, p}(\mu)$; and

(ii) $Z = (1 + {\mathcal O}(\varepsilon)) {\mathbb X} + {\mathcal
O}_{\theta, p}(1) + {\mathcal O}_{{\mathbb X}, \theta, p}(\mu)$.
\end{prop}
\noindent {\bf Proof:} \ The proof of (a) is identical to that of
Lemma 3.1 and 3.3 of \cite{WOk} and the proof of (b) is identical
to that of Lemma 3.4 of \cite{WOk}. \hfill $\square$

\medskip

\noindent {\bf B. The induced map ${\mathcal N}: \Sigma^+ \to
\Sigma^-$.} \ For $({\mathbb X}, {\mathbb Y}, \theta) \in
\Sigma^+$ we have ${\mathbb X} = \varepsilon \mu^{-1}$ by
definition. Similarly, for $({\mathbb X}, {\mathbb Y}, \theta) \in
\Sigma^-$ we have ${\mathbb Y}= \varepsilon \mu^{-1}$.  Denote a
point on $\Sigma^+$ by using $({\mathbb Y}, \theta)$ and a point
on $\Sigma^-$ by using $({\mathbb X}, \theta)$, and let
$$
({\mathbb X}_1, \theta_1) = {\mathcal N}({\mathbb Y}, \theta)
$$
for $({\mathbb Y}, \theta)\in \Sigma^+$.
\begin{prop}\label{prop1-s2.2b}
We have for $({\mathbb Y}, \theta) \in \Sigma^+$,
\begin{equation}\label{f1-s2.2b}
\begin{split}
{\mathbb X}_1 & =  (\mu \varepsilon^{-1})^{\frac{\alpha}{\beta}-1}
{\mathbb Y}^{\frac{\alpha}{\beta}} \\
\theta_1 & = \theta + \frac{\omega} {\beta} \ln (\varepsilon
\mu^{-1}) - \frac{\omega}{\beta} \ln {\mathbb Y}.
\end{split}
\end{equation}
\end{prop}
\noindent {\bf Proof:} \ Let $T$ be the time it takes for the
solution of the linearized equation of Proposition
\ref{prop1-s2.1b} from $(\varepsilon, Y, \theta) \in \Sigma^+$ to
get to $(X_1, \varepsilon, \theta_1) \in \Sigma^-$. We have
\begin{equation*}
X_1 = \varepsilon e^{-\alpha T}, \ \ \ \ \varepsilon = Y e^{\beta
T}, \ \ \ \ \theta_1 = \theta + \omega T,
\end{equation*}
from which (\ref{f1-s2.2b}) follows. \hfill $\square$

\medskip

\noindent {\bf C. The induced map ${\mathcal M}: \Sigma^- \to
\Sigma^+$.} \ Let ${\mathbb H}_i(s)$, $i = 1, 2$ be as in
(\ref{f4-s2.1c}). In what follows, we write
\begin{equation}\label{f1-s2.2c}
\begin{split}
A_L & = \int_{-L^-}^{L^+} {\mathbb H}_1(s) e^{- \int_{0}^{s}
E(\tau)
d\tau} ds \\
\phi_L(\theta) & = \int_{-L^-}^{L^+} {\mathbb H}_2(s) \sin(\theta
+ \omega s + \omega L^-) e^{- \int_{0}^{s} E(\tau) d\tau} ds.
\end{split}
\end{equation}
We also write
\begin{equation} \label{f2-s2.2c}
P_L  = e^{\int_{-L^-}^{L^+} E(s) ds}, \ \ \ P_L^+ =
e^{\int_{0}^{L^+} E(s) ds}.
\end{equation}
Note that for $P_L$ we integrate from $s = -L^-$ to $s = L^+$,
while for $P_L^+$ the integration starts from $s = 0$. First we
have \
\begin{lemma}\label{lem1-s2.2c}
$$
P_L \sim \varepsilon^{\frac{\alpha}{\beta} - \frac{\beta}{\alpha}}
<< 1, \ \ \ \ \ \ \ P_L^+ \sim \varepsilon^{-\frac{\beta}{\alpha}}
>> 1.
$$
\end{lemma}
\noindent {\bf Proof:} \ Same as Lemma 3.5 in \cite{WOk}. \hfill
$\square$

\smallskip

For $q = (s^-, Z, \theta) \in \Sigma^-$, the value of $s^-$ is
uniquely determined by that of $(Z, \theta)$ through Proposition
\ref{prop1-s2.2b}(b)(i). So we can use $(Z, \theta)$ to represent
$q$. Let $(s(t), Z(t), \theta(t))$ be the solution of equation
(\ref{f2-s2.1c}) initiated at $(s^-, Z, \theta)$, and $\tilde t$
be the time $(s(\tilde t), Z(\tilde t), \theta(\tilde t))$ hit
$\Sigma^+$. By definition ${\mathcal M}(q) = (s(\tilde t),
Z(\tilde t), \theta(\tilde t))$. In what follows we write
$$
\hat Z = Z(\tilde t), \ \  \hat \theta = \theta(\tilde t).
$$
\begin{prop}\label{prop1-s2.2c}
Denote $(\hat Z, \hat \theta) = {\mathcal M}(Z, \theta)$. We have
\begin{equation}\label{f1-s2.2c}
\begin{split}
\hat Z & =  P_L^+ (\rho A_L + \phi_L(\theta))+ P_L Z + {\mathcal
O}_{Z, \theta, p}(\mu) \\
\hat \theta & = \theta + \omega (L^+ + L^-) + {\mathcal O}_{Z,
\theta, p}(\mu) .
\end{split}
\end{equation}
\end{prop}
\noindent {\bf Proof:} \ The same as that of Proposition 3.2 of
\cite{WOk}. \hfill $\square$

\medskip

\noindent {\bf D. The return map ${\mathcal F} = {\mathcal N}
\circ {\mathcal M}$.} \ We are now ready to compute the return map
${\mathcal F} = {\mathcal N} \circ {\mathcal M}: \Sigma^- \to
\Sigma^-$. We use $({\mathbb X}, \theta)$ to represent a point on
$\Sigma^-$ and denote $({\mathbb X}_1, \theta_1) = {\mathcal
F}({\mathbb X}, \theta)$.
\begin{prop}\label{prop1-s2.2d}
The map ${\mathcal F} ={\mathcal N} \circ {\mathcal M}: \Sigma^-
\to \Sigma^-$ is given by
\begin{equation}\label{f1-s2.2d}
\begin{split}
{\mathbb X}_1 & =  (\mu \varepsilon^{-1})^{\frac{\alpha}{\beta}-1}
[(1 + {\mathcal O}(\varepsilon))P_L^+ {\mathbb F}({\mathbb X},
\theta)]^{\frac{\alpha}{\beta}} \\
\theta_1 & = \theta + \omega (L^+  + L^-) + \frac{\omega}{\beta}
\ln \mu^{-1} \varepsilon (1 + {\mathcal O}(\varepsilon))P_L^+ -
\frac{\omega}{\beta} \ln {\mathbb F}({\mathbb X}, \theta) +
{\mathcal O}_{{\mathbb X}, \theta, p}(\mu)
\end{split}
\end{equation}
where
\begin{equation}\label{f2-s2.2d}
\begin{split}
{\mathbb F}({\mathbb X}, \theta) & = (\rho A_L + \phi_L(\theta))+
P_L (P_L^+)^{-1} (1+ {\mathcal O}(\varepsilon)){\mathbb X} \\
&  \ \ \ \ + (P_L^+)^{-1}(1+P_L) {\mathcal O}_{\theta, p}(1) +
{\mathcal O}_{{\mathbb X}, \theta, p}(\mu),
\end{split}
\end{equation}
and $P_L, P_L^+$ and $\phi_L(\theta)$ are as in (\ref{f1-s2.2c})
and (\ref{f2-s2.2c}).
\end{prop}
\noindent {\bf Proof:} \ By using Proposition \ref{prop1-s2.2c}
and Proposition \ref{prop1-s2.2a}(b), we have
\begin{equation*}
\begin{split}
\hat Z & = P_L(1+ {\mathcal O}(\varepsilon)){\mathbb X} + P_L^+
(\rho A_L + \phi_L(\theta)) + P_L {\mathcal O}_{\theta, p}(1) +
{\mathcal O}_{{\mathbb X}, \theta, p}(\mu) \\
\hat \theta & = \theta + \omega (L^+ + L^-) + {\mathcal
O}_{{\mathbb X}, \theta, p}(\mu).
\end{split}
\end{equation*}

Let $\hat {\mathbb Y}$ be the ${\mathbb Y}$-coordinate for $(\hat
Z, \hat \theta)$, we have from Proposition \ref{prop1-s2.2a}(a),
$$
\hat {\mathbb Y} = (1 + {\mathcal O}(\varepsilon))P_L^+ {\mathbb
F}({\mathbb X}, \theta)
$$
where ${\mathbb F}({\mathbb X}, \theta)$ is as in
(\ref{f2-s2.2d}). We then obtain (\ref{f1-s2.2d}) by using
(\ref{f1-s2.2b}).  \hfill $\square$

\section{Proofs of Theorems \ref{th1}-\ref{th4}}\label{s3}
Recall that $\ell = \ell_{\lambda}$ is the homoclinic solution of
equation of Proposition \ref{prop1} for equation (\ref{f1-s1}),
which we write as $\ell(s) = (a(s), b(s))$. $(u(s), v(s))$ is the
unit tangent vector of $\ell$ at $\ell(s)$. Also recall that
${\mathbb H}_i(s)$, $i = 1, 2$ are as in (\ref{f4-s2.1c}) and
$E(s)$ is as in (\ref{f3-s2.1c}). Let
\begin{equation}\label{f1-s3}
\begin{split}
A & = \int_{-\infty}^{\infty}{\mathbb H}_1(s) e^{- \int_{0}^{s}
E(\tau) d\tau} ds  \\
C(\omega) & = \int_{-\infty}^{\infty} {\mathbb H}_2(s)
\cos (\omega s) e^{- \int_{0}^{s} E(\tau) d\tau} ds \\
S(\omega) & = \int_{-\infty}^{\infty}{\mathbb H}_2(s) \sin (\omega
s) e^{- \int_{0}^{s} E(\tau) d\tau} ds.
\end{split}
\end{equation}
First we need to estimate $A, C(\omega)$ and $S(\omega)$.
\begin{prop}\label{prop1-s3}
There exists a $\lambda_0(R_{\omega})$, sufficiently small
depending on $R_{\omega}$, such that for all $\omega \in (0,
R_{\omega})$ and all $\lambda \in (0, \lambda_0)$, we have

(i) $A > 1$; and

(ii) $\sqrt{C^2(\omega) + S^2(\omega)} > \left(e^{-\frac{1}{2}
\omega \pi} + e^{\frac{1}{2} \omega \pi}\right)^{-1}$.
\end{prop}
Proposition \ref{prop1-s3} is proved in Section \ref{s4}.

\medskip

\noindent {\bf The return maps:} \ To apply Proposition
\ref{prop1-s2.2d}, we first fix a $R_{\omega} >> 100 \beta$. Then
we pick a $\lambda \in (0, \lambda_0(R_{\omega}))$ satisfying the
non-resonance assumption (H). We then chose $R_{\rho},
\varepsilon, R_{\mu}$ such that
$$
R_{\mu} >> \varepsilon^{-1} >> R_{\rho} >> \lambda^{-1}.
$$

Recall that $p = (\omega, \rho, \ln \mu) \in {\mathbb P}$
represents the forcing parameter and
$$
\Sigma^-  = \{ (\theta, {\mathbb X}): \ \ \theta \in {\mathbb
R}/(2 \pi {\mathbb Z}), \ \ |{\mathbb X}| < 1 \}.
$$
Let $(\theta_1, {\mathbb X}_1) = {\mathcal F}(\theta, {\mathbb
X})$ for $(\theta, {\mathbb X}) \in \Sigma^-$ where ${\mathcal F}$
is from Proposition \ref{prop1-s2.2d}. We have
\begin{equation}\label{f1-s3.1a}
\begin{split}
\theta_1 & = \theta + {\bf a} - \frac{\omega}{\beta} \ln {\mathbb
F}(\theta, {\mathbb X}, p) \\
{\mathbb X}_1 & =  {\bf b} [{\mathbb F}(\theta, {\mathbb X},
p)]^{\frac{\alpha}{\beta}}
\end{split}
\end{equation}
where
\begin{equation} \label{f2-s3.1a}
\begin{split}
{\bf a} & = \frac{\omega}{\beta} \ln \mu^{-1} + \omega (L^+  +
L^-) + \frac{\omega}{\beta} \ln(\varepsilon (1 +
{\mathcal O}(\varepsilon))P_L^+ A_L) \\
{\bf b} & = (\mu \varepsilon^{-1})^{\frac{\alpha}{\beta}-1} [(1 +
{\mathcal O}(\varepsilon))P_L^+ A_L]^{\frac{\alpha}{\beta}}  \\
\end{split}
\end{equation}
and
\begin{equation} \label{f3-s3.1a}
{\mathbb F}(\theta, {\mathbb X}, p) = \rho + {\bf c} \sin \theta +
{\bf k} {\mathbb X} + {\mathbb E}(\theta, p) + {\mathcal
O}_{\theta, {\mathbb X}, p}(\mu),
\end{equation}
in which
\begin{equation} \label{f4-s3.1a}
\begin{split}
{\bf c} & = (A_L)^{-1} \sqrt{C_L^2 + S_L^2} \\
{\bf k} & =  P_L (A_L P_L^+)^{-1} (1+ {\mathcal O}(\varepsilon))
\end{split}
\end{equation}
and
\begin{equation} \label{f5-s3.1a}
{\mathbb E}(\theta, p)  = (A_L P_L^+)^{-1}(1+P_L) {\mathcal
O}_{\theta, p}(1).
\end{equation}
Note that in getting (\ref{f1-s3.1a}) we have changed $\theta +
\omega L^- + c_0$ to $\theta$ where $c_0$ is such that $\tan c_0 =
C_L^{-1} S_L$. ${\bf a}, {\bf b}, {\bf c}, {\bf k}$ and ${\mathbb
E}(\theta, p)$ are as follows:

\smallskip

(1) As parameters, ${\bf b}$ and ${\bf a}$ are functions of $p=
(\omega, \rho, \ln \mu)$ and $\varepsilon$. If we regard $\omega,
\rho$ and $\varepsilon$ as been fixed, and ${\mathcal F}_{\mu} =
{\mathcal F}$ as a one parameter family, then ${\bf b} \to 0$ as
$\mu \to 0$. We can then think ${\mathcal F}$ as an unfolding of
the 1D maps
$$
f(\theta) = \theta + {\bf a} - \frac{\omega}{\beta} \ln (\rho +
{\bf c} \sin \theta + {\mathbb E}(\theta, 0))
$$
where ${\bf a} \approx - \omega \beta^{-1} \ln \mu$. ${\bf a} \to
+\infty$ as $\mu \to 0$.

\smallskip

(2) ${\bf a}$ is a large number. But since it appears in the
angular component we can module it by $2 \pi$. Varying $\mu$ from
$R_{\mu}^{-1}>0$ to zero is to run ${\bf a}$ over $({\bf a}_0,
+\infty)$ for some ${\bf a}_0 \sim \omega \beta^{-1} \ln R_{\mu}$.

\smallskip

(3) ${\bf c}$, ${\bf k}$ are independent of forcing parameters $p$
but are functions of $\varepsilon$. From Proposition
\ref{prop1-s3} and Lemma \ref{lem1-s2.2c} we have
$$
{\bf c} \neq 0, \ \ \ \ \ {\bf k} \sim
\varepsilon^{\frac{\alpha}{\beta}}
$$
provided that $\varepsilon$ is sufficiently small. Recall that the
value of $\varepsilon$ is selected after that of $R_{\omega}$ and
$R_{\rho}$.

\smallskip

(4) From Lemma \ref{lem1-s2.2c} and Proposition \ref{prop1-s3} we
have
$$
{\mathbb E}(\theta, \mu) \sim \varepsilon^{\beta \alpha^{-1}}
{\mathcal O}_{\theta, p}(1).
$$
When $\varepsilon$ is sufficiently small, ${\mathbb E}(\theta,
\mu)$ is a $C^r$-small perturbation to $\rho + {\bf c} \sin
\theta$.

\medskip

\noindent {\bf Proof of Theorem \ref{th1}:} \ To prove theorem
\ref{th1}, we fix the values of $\omega, \mu$ and $\varepsilon$
and regard $\rho \in (0, R_{\rho})$ as the parameter in ${\mathcal
F}$. Let
$$
F(\theta, \rho) = {\mathbb F}(\theta, 0, p) = \rho + {\bf c} \sin
\theta + {\mathbb E}(\theta, p) + {\mathcal O}_{\theta, p}(\mu)
$$
and denote
$$
M(\rho) = \min_{\theta \in S^1} F(\theta, \rho).
$$
$M(\rho)$ is a continuous function of $\rho$.

We prove that $M(\rho)$ is monotonic in $\rho$. For this it
suffices to prove that for all fixed $\theta \in S^1$,
$$
\partial_{\rho} F(\theta, \rho) = 1 + \partial_{\rho} E(\theta, p)
+ \partial_{\rho} {\mathcal O}_{\theta, p}(\mu) >  1 - K
\varepsilon^{\frac{\alpha}{\beta}} - K(\varepsilon) \mu > 0.
$$

Next we prove that $W^s \cap W^u = \emptyset$ if and only if
$M(\rho) > 0$. Let $l^u_{loc}$ be the curve in $\Sigma^-$ defined
by ${\mathbb X} = 0$. From (\ref{f1-s3.1a}) we know that
${\mathcal F}^n(l_{loc}^u)$ is well-defined for all $n > 0$ if and
only if $M(\rho) > 0$. We also know that $W^u \cap W^u =
\emptyset$ if and only if ${\mathcal F}^n(l_{loc}^u)$ is
well-defined for all $n > 0$ because $W^u$ came out of
$l_{loc}^u$.

To finish the proof of Theorem \ref{th1}, it now suffices to
observe that $M(R_{\rho}) > 0$ and $M(0) < 0$ because ${\bf c}
\neq 0$. The surface $S^*$ in ${\mathbb P}$ is defined by $M(\rho)
= 0$. In fact we have
\begin{equation}\label{S-star}
S^*(\omega, \mu) = {\bf c} + {\mathcal O}(\varepsilon + \mu)
\end{equation}
from $M(\rho) = 0$. \hfill $\square$

\medskip

\noindent {\bf Proof of Theorem \ref{th2}:} \ The formulas
(\ref{f1-s3.1a}) for the return maps are identical to these
obtained in Sect. 4.1 of \cite{WOk}, and Proposition
\ref{prop1-s3} implies that assumption (H2) in Sect. 4.1 of
\cite{WOk} holds.\footnote{Though there is a difference in
writing: in this paper we kept a factor $\rho$ inside of ${\mathbb
F}$ while in \cite{WOk} we gave it to ${\bf b}$ and ${\bf a}$.
With the current writing the proof of Theorem \ref{th1} is easier
to present.} Assumption (H1) in Sect. 2.1 of \cite{WOk} is (H).
Consequently, Theorems 1-3 of \cite{WOk} apply. A specific choice
of the domain $DH$ is defined in the paragraph entitled
``specifications of parameters" of Sect. 4.1 of \cite{WOk}.
Theorem \ref{th2}(i) is Theorem 1 of \cite{WOk}; Theorem
\ref{th2}(ii) is Theorem 2 of \cite{WOk}; and Theorem
\ref{th2}(iii) is Theorem 3 of \cite{WOk}. \hfill $\square$

\medskip

\noindent {\bf Proof of Theorem \ref{th3}:} \ Let $\rho_0 > 0$ be
such that
$$
\rho_0 > 2 A^{-1} \sqrt{C^2(0) + S^2(0)}
$$
where $C(0), S(0)$ are the values of $C(\omega)$ and $S(\omega)$
at $\omega = 0$. For $(\omega, \rho, \mu) \in {\mathbb P}$
satisfying $\rho
> \rho_0$, let
$$
M(\omega, \rho, \mu) = \min_{\theta \in S^1} {\mathbb F}(\theta,
0, p).
$$
We have
\begin{equation}\label{f6-s3.1a}
M(0, \rho, \mu) > \frac{1}{5} \rho_0 > 0.
\end{equation}
We now define $Q$ by letting
$$
Q(\rho, \mu) = \min\{ \omega>0: \ \ \ \omega \leq 10^{-5}
M(\omega, \rho, \mu) {\bf c}^{-1}(\omega); \ \ M(\omega, \rho,
\mu) \geq \frac{1}{10} \rho_0 \}.
$$
$Q(\rho, \mu)$ is well defined because of (\ref{f6-s3.1a}) and the
fact that $M(\omega, \rho, \mu)$ is continuous in $\omega$. In the
rest of this proof we assume $(\omega, \rho, \mu) \in {\mathbb P}$
is such that $\omega < Q(\rho, \mu)$.

With (\ref{f1-s3.1a}) for ${\mathcal F}$ and the way $Q(\rho,
\mu)$ is defined above, our proof of Theorem \ref{th3} is the same
as that of Theorem 1 of \cite{WY4}. See Sect. 4.1 of \cite{WY4}.
For the existence of an attracting invariant curve, we observe
that $D{\mathcal F}$ has a dominating splitting on $\Sigma^-$
(Sect. 4.1.1 of \cite{WY4}). Theorem \ref{th3}(2) then follows
from Denjoy theory and a direct application of a Theorem of Herman
\cite{He} to the circle diffeomorphisms induced on the attracting
invariant curves (Sect. 4.1.2 of \cite{WY4}).

To be more precise, we let ${\bf v} = (u, v)$ be a tangent vector
of $\Sigma^-$ at $q$ and let $s({\bf v}) = v u^{-1}$ for $q =
(\theta, {\mathbb X}) \in \Sigma^-$. $s({\bf v})$ is the slope of
${\bf v}$. Let ${\mathcal C}_h(q)$ be the collection of all ${\bf
v}$ satisfying $|s({\bf v})|< \frac{1}{100}$, and ${\mathcal
C}_v(q)$ be the collection of all ${\bf v}$ satisfying $|s({\bf
v})| > 100$. We first prove that
\begin{itemize}
\item[(i)] on $\Sigma^-$, we have $D{\mathcal F}_q({\mathcal
C}_h(q)) \subset {\mathcal C}_h({\mathcal F}(q))$, and
$|D{\mathcal F}_q({\bf v})| > \frac{1}{2}|{\bf v}|$; \item[(ii)]
on ${\mathcal F}(\Sigma^-)$, we have $D{\mathcal
F}^{-1}_q({\mathcal C}_v(q)) \subset {\mathcal C}_v({\mathcal
F}(q))$ and $|D{\mathcal F}^{-1}_q({\bf v})| > 100 |{\bf v}|$ for
all ${\bf v} \in {\mathcal C}_h(q)$
\end{itemize}

To prove (i) we compute $D{\mathcal F}$ by using (\ref{f1-s3.1a}).
Let $(\theta_1, {\mathbb X}_1) = {\mathcal F}(\theta, {\mathbb
X})$, we have
\begin{equation}\label{f1-add-s3.2}
D{\mathcal F} = \left( \begin{array}{cc} \frac{\partial
\theta_1}{\partial \theta} & \frac{\partial \theta_1}{\partial
{\mathbb X}}
\\
\frac{\partial {\mathbb X}_1}{\partial \theta} & \frac{\partial
{\mathbb X}_1}{\partial {\mathbb X}}
\end{array} \right) = \left(
\begin{array}{cc} 1
- \omega \beta^{-1} \frac{1}{{\mathbb F}} \frac{\partial {\mathbb
F}}{\partial \theta} &  \omega \beta^{-1} \frac{1}{{\mathbb F}}
\frac{\partial {\mathbb F}}{\partial {\mathbb X}} \\
\alpha \beta^{-1}  {\bf b} {\mathbb F}^{\alpha \beta^{-1}-1}
\frac{\partial {\mathbb F}}{\partial \theta} & \alpha \beta^{-1}
{\bf b} {\mathbb F}^{\alpha \beta^{-1} -1} \frac{\partial {\mathbb
F}}{\partial {\mathbb X}}
\end{array} \right)
\end{equation}
where ${\mathbb F} = {\mathbb F}(\theta, {\mathbb X}, p)$ is as in
(\ref{f3-s3.1a}) and
\begin{equation*}
\begin{split}
\frac{\partial {\mathbb F}}{\partial \theta} & = {\bf c} \cos
\theta + \varepsilon^{\beta \alpha^{-1}} {\mathcal O}_{\theta,
p}(1)
+ {\mathcal O}_{\theta, {\mathbb X}, {\bf a}}(\mu)\\
\frac{\partial {\mathbb F}}{\partial {\mathbb X}} & = {\bf k} +
{\mathcal O}_{\theta, {\mathbb X}, p}(\mu).
\end{split}
\end{equation*}
(i) then follows from ${\bf b} \approx \mu^{\lambda \beta^{-1}} <<
1$ and
\begin{equation}\label{f7-s3.1a}
\left| \omega \beta^{-1} \frac{1}{\mathbb F} \frac{\partial
{\mathbb F}}{\partial \theta} \right|< 10^{-4}
\end{equation}
where (\ref{f7-s3.1a}) is from $\omega < Q(\rho, \mu)$. To prove
(ii) we use again (\ref{f7-s3.1a}) and
\begin{equation}\label{f4-add-s3.2}
D{\mathcal F}^{-1} = \frac{1}{\alpha \beta^{-1} {\bf b} {\mathbb
F}^{\alpha \beta^{-1} -1} \frac{\partial {\mathbb F}}{\partial
{\mathbb X}}} \left(
\begin{array}{cc} \alpha \beta^{-1}
{\bf b} {\mathbb F}^{\alpha \beta^{-1} -1} \frac{\partial {\mathbb
F}}{\partial {\mathbb X}} &  - \omega \beta^{-1} \frac{1}{{\mathbb
F}}
\frac{\partial {\mathbb F}}{\partial {\mathbb X}} \\
- \alpha \beta^{-1}  {\bf b} {\mathbb F}^{\alpha \beta^{-1}-1}
\frac{\partial {\mathbb F}}{\partial \theta} & 1 - \omega
\beta^{-1} \frac{1}{{\mathbb F}} \frac{\partial {\mathbb
F}}{\partial \theta}
\end{array} \right).
\end{equation}
With both (i) and (ii), Theorem \ref{th3}(1) follows now from the
standard graph transformation argument. Proof for Theorem
\ref{th3}(2) is the same as the proof of Theorem 1(a)(b) of
\cite{WY4}. See Sect. 4.1.2 of \cite{WY4}, where $T$ is the
parameter ${\bf a}$ for here. \hfill $\square$

\medskip

\noindent {\bf Proof of Theorem \ref{th4}:} \ Theorem \ref{th4}(2)
is proved by applying the theory of rank one maps of \cite{WY1}
and \cite{WY2} to ${\mathcal F}$ of (\ref{f1-s3.1a}). The domain
$DC$ is such that
$$
\frac{202}{99} \frac{\sqrt{C^2(\omega) + S^2(\omega)}}{A} < \rho <
\frac{396}{101} \frac{\sqrt{C^2(\omega) + S^2(\omega)}}{A}.
$$
We also need to assume that $\omega$ is sufficiently large for
$(\omega, \rho) \in DC$ in order to prove that ${\mathcal
F}_{\mu}$ is an admissible family. For more details see Section 6
of \cite{WO}.

Theorem \ref{th4}(1) is a specific case of Theorem 3 of \cite{LW}.
\hfill $\square$

\section{Proof of Proposition \ref{prop1-s3}} \label{s4}
In this section we prove Proposition \ref{prop1-s3}. We start with
the conservative part of the autonomous equation (\ref{f1-s1}),
that is, the equation
\begin{equation}\label{f1-s4}
\frac{d^2 q}{dt^2} - q + q^3 = 0.
\end{equation}

Denote $p = \frac{d q}{dt}$. We write equation (\ref{f1-s4})as
\begin{equation}\label{f2-s4}
\frac{d q}{dt} = p, \ \ \ \ \ \ \frac{d p}{ dt} = q - q^3.
\end{equation}
With a simple linear change of coordinates
\begin{equation*}
x  =  \frac{1}{2} (q -  p),  \ \ \ \  y  = \frac{1}{2} (q + p),
\end{equation*}
we write equation (\ref{f1-s4}) in $(x, y)$ as
\begin{equation}\label{f3-s4}
\frac{dx}{dt} = - x + \frac{1}{2} (x+y)^3, \ \ \ \ \ \ \
\frac{dy}{dt} = y - \frac{1}{2} (x+y)^3.
\end{equation}
Let
\begin{equation*}
a(t)  = \frac{ 2\sqrt{2} e^{3t}}{(1 + e^{2t})^2}, \ \ \ \ \ \ \
b(t) =\frac{2\sqrt{2}e^{t}}{(1 + e^{2t})^2}.
\end{equation*}
$(x, y) = (0, 0)$ is a homoclinic saddle and $(x, y) = (a(t),
b(t))$ is a homoclinic solution of equation (\ref{f3-s4})
initiated at $(\sqrt{2}, \sqrt{2})$. Denote
$$
\ell = \{ \ell(t) = (a(t), b(t)), \ t \in {\mathbb R} \}.
$$

It follows that
\begin{equation}\label{f4-s4}
\begin{split}
u(t) &= \frac{-(e^{2 t} - 3) }{\sqrt{(e^{2t} - 3)^2
+ (e^{-2 t} - 3)^2}} \\
v(t) &= \frac{e^{-2 t} - 3}{\sqrt{(e^{2 t} - 3)^2 + (e^{-2 t} -
3)^2}},
\end{split}
\end{equation}
where
$$
(u(t), v(t)) = \left|\frac{d}{dt} \ell(t)\right|^{-1} \frac{d}{dt}
\ell(t)
$$
is the unit tangent vector of $\ell$ at $\ell(t)$.

Let us re-write equation (\ref{f3-s4}) as
\begin{equation}\label{f5-s4}
\frac{dx}{dt} = - x + f(x, y) \ \ \ \ \ \ \ \ \frac{dy}{dt} = y +
g(x, y)
\end{equation}
where
$$
f(x,y) = \frac{1}{2} (x+y)^3, \ \ \ \ \ \ g(x, y) = - \frac{1}{2}
(x+y)^3.
$$
Let
\begin{equation*}
\begin{split}
E(t) & = v^2(t) (- 1 + \partial_x f(a(t), b(t))) + u^2(t)(1 +
\partial_y g(a(t), b(t))) \\
& \ \ - u(t) v(t) (\partial_y f(a(t), b(t)) + \partial_x g(a(t),
b(t))).
\end{split}
\end{equation*}
We have
\begin{equation}\label{f6-s4}
E(t) = - \frac{(e^{-2 t} - 3)^2 - (e^{2 t} - 3)^2}{(e^{2 t} - 3)^2
+ (e^{-2  t} - 3)^2}\left( 1 - \frac{12 e^{ 2 t}}{(1+e^{2 t})^2}
\right).
\end{equation}

Let
\begin{equation}\label{f5-s2.1}
K(s) = - \int_0^s E(s) ds.
\end{equation}

\begin{lemma}\label{lem1-s4}
For $s \in (-\infty, \infty)$,
$$
K(s) = \frac{1}{2} \ln \frac{8 e^{2 s} ((1 - 3e^{2 s})^2 + e^{4 s}
(e^{2 s}- 3)^2)}{(e^{2s} + 1)^6}.
$$
\end{lemma}
\noindent {\bf Proof:} \ Let $x = e^{2t}$. We have\footnote{The
author would like to thank Ali Oksasoglu for making this
computation.}
$$
K(s)  =  \frac{1}{2} \int_{1}^{e^{2s}} \frac{1}{x}\frac{(1 - 3x)^2
- x^2(x - 3)^2}{(1 - 3x)^2 + x^2(x- 3)^2} \left(1 - \frac{12 x}{(1
+ x)^2}\right) d x.
$$
Observe that
\begin{eqnarray*}
(1 - 3x)^2 + x^2(x- 3)^2  = ((x-a)^2 + a^2)((x-b)^2 + b^2)
\end{eqnarray*}
where $a, b >0$ satisfying $a + b = 3, \ ab = \frac{1}{2}$. We
have
\begin{eqnarray*}
K(s) & = & \frac{1}{2} \int_1^{e^{2s}} \frac{(-x^4 + 6x^3 - 6x
+1)(x^2 - 10 x +1)}{x(x+1)^2 ((x-a)^2 + a^2)((x-b)^2 + b^2) } dx
\\
& = & \frac{1}{2} \int_1^{e^{2s}} \left(\frac{1}{x} -
\frac{6}{1+x} + \frac{4x^3 - 18 x^2 + 36 x -6}{((x-a)^2 + a^2)((x-b)^2 + b^2)} \right) dx \\
& = & \frac{1}{2} \int_1^{e^{2s}} \left(\frac{1}{x} -
\frac{6}{1+x} + \frac{2(x-a)}{(x-a)^2 + a^2} +
\frac{2(x-b)}{(x-b)^2 + b^2} \right) dx.
\end{eqnarray*}
We then have
\begin{eqnarray*}
K(s) & = &\frac{1}{2} \ln \frac{8 e^{2s} ((e^{2s}-a)^2 +
a^2)((e^{2s}-b)^2 + b^2)}{(e^{2s} + 1)^6} \\
& = & \frac{1}{2} \ln \frac{8 e^{2s} ((1 - 3e^{2s})^2 +
e^{4s}(e^{2s}- 3)^2)}{(e^{2s} + 1)^6}.
\end{eqnarray*}

\hfill $\square$

\smallskip

Through direct evaluation using Lemma \ref{lem1-s4}, we obtain the
values of the following definite integrals as
\begin{equation}\label{f7-s4}
\int_{-\infty}^{\infty} (u(s)+v(s))(b(s) + (a(s))^2(b(s)-a(s))
e^{-\int_0^s E(t) dt} ds  = \frac{16}{15}.
\end{equation}

\smallskip

We have two more integrals to evaluate, and they are
\begin{equation}\label{f8-s4}
\begin{split}
C(\omega) & = \int_{-\infty}^{\infty} (u(s)+v(s))(b(s) + a(s))^2
\cos (\omega s) e^{-\int_{0}^s E(t) dt} ds,  \\
S(\omega) & = \int_{-\infty}^{\infty} (u(s)+v(s))(b(s) + a(s))^2
\sin (\omega s) e^{-\int_{0}^s E(t) dt} ds.
\end{split}
\end{equation}

\begin{lemma}\label{lem2-s4} We have
\begin{equation*}
C(\omega) =  \frac{16 \sqrt{2} \pi  }{e^{-\frac{1}{2} \omega \pi}
+ e^{\frac{1}{2} \omega \pi}}, \ \ \ \ \  S(\omega)= -
\frac{2\sqrt{2}}{3} \cdot \frac{ \omega(1 +
\omega^2)}{e^{-\frac{1}{2} \omega \pi} + e^{\frac{1}{2} \omega
\pi}}.
\end{equation*}
\end{lemma}
\noindent {\bf Proof:} \ Using Lemma \ref{lem1-s4}, we have
$$
C(\omega)  = 16\sqrt{2} I_c,  \ \ \ \ \ \ \ S(\omega)  = 16
\sqrt{2} I_s
$$
where
\begin{equation*}
\begin{split}
I_c & = \int_{-\infty}^{\infty} \frac{e^{3 s}(1- e^{2 s})}{(e^{2
s} +1)^4} \cdot \cos(\omega
s) \cdot ds \\
I_s & = \int_{- \infty}^{\infty} \frac{e^{3 s}(1- e^{2 s})}{(e^{2
s} +1)^4}  \cdot \sin (\omega s) \cdot ds.
\end{split}
\end{equation*}
Let $I = I_c  + i I_s$.
$$
I = \int_{-\infty}^{\infty} \frac{e^{3 s}(1-e^{2 s})}{(e^{2 s}
+1)^4} \cdot e^{i \omega s} \cdot ds.
$$

We evaluate $I$ as follows. On the $z = x + i y$ plan, let
$$
\ell_1  =  \{x,  \ \ x: - \infty \to \infty \},  \ \ \ \ \ell_2  =
\{x + \pi i, \ \ x: \infty  \to  - \infty \}
$$
and
$$
f(z) = \frac{e^{3 z}(1-e^{2 z})}{(e^{2 z} +1)^4} \cdot e^{i \omega
z}.
$$
We have
$$
\int_{\ell_1} f(z) dz = I, \ \ \ \ \ \ \ \int_{\ell_2} f(z) d z =
e^{-\omega \pi} I.
$$
So by the residue theorem,
$$
(1+e^{-\omega \pi})I = 2 \pi i Res(f(z))_{z = \frac{\pi i}{2}}.
$$
Let $t = z - \frac{\pi i}{2}$, we write $f(z)$ as
\begin{eqnarray*}
f(z) & = & - i e^{-\frac{1}{2} \omega \pi} \frac{e^{3 t}(e^{2 t} +
1)}{(1- e^{2
t})^4} \cdot e^{i \omega t} \\
& = & - \frac{i}{16} e^{-\frac{1}{2} \omega \pi} \frac{e^{(5+
\omega i) t} + e^{(3 + \omega i)t}}{t^4(1 + t + \frac{2}{3} t^2 +
\frac{1}{3} t^3 + {\mathcal
O}(t^4))^4}.  \\
\end{eqnarray*}

We have
\begin{eqnarray*}
e^{(5+ i \omega)t} + e^{(3+ i \omega)t} & = & 1 + (8 + 2 i
\omega)t  + (17 + 8 i \omega - \omega^2) t^2 \\
& & \ \  + \frac{1}{3}(76 - 12 \omega^2 + (51 \omega - \omega^3)i)
t^3 + {\mathcal O}(t^4).
\end{eqnarray*}
We also have
$$
\frac{1}{(1+x)^4} = 1- 4x + 10 x^2 - 20 x^3 + \cdots
$$
so
\begin{eqnarray*}
(1 + t + \frac{2}{3} t^2 + \frac{1}{3} t^3)^{-4} & = & 1 - 4 t
+\frac{22}{3} t^2 - 8 t^3 +{\mathcal O}(t^4).
\end{eqnarray*}
Consequently,
$$
Res(f(z)) = - \frac{i}{16} e^{-\frac{1}{2} \omega \pi}(8 -
\frac{1}{3}\omega (1 + \omega^2) i),
$$
from which it follows that
$$
I = \frac{1}{8} \frac{\pi e^{-\frac{1}{2} \omega \pi} }{1 +
e^{-\omega \pi}}(8 - \frac{1}{3}\omega (1 + \omega^2) i).
$$
In conclusion, we have
\begin{equation*}
C(\omega) =  \frac{16 \sqrt{2} \pi  }{e^{-\frac{1}{2} \omega \pi}
+ e^{\frac{1}{2} \omega \pi}}, \ \ \ \ \  S(\omega)= -
\frac{2\sqrt{2}}{3} \cdot \frac{ \omega(1 +
\omega^2)}{e^{-\frac{1}{2} \omega \pi} + e^{\frac{1}{2} \omega
\pi}}.
\end{equation*}

\hfill $\square$

\medskip

\noindent {\bf Proof of Proposition \ref{prop1-s3}} \ We use
subscript $\lambda$ to stress the role of $\lambda$. We denote
equation (\ref{f4-s2.1a}) as (\ref{f4-s2.1a})$_{\lambda}$, and the
homoclinic solution $\ell_{\lambda}$ of Proposition \ref{prop1} as
$$
x = a_{\lambda}(t), \ \ \ \ \ \ y = b_{\lambda}(t).
$$
Similarly, $E_{\lambda}$ is defined through (\ref{f3-s2.1c}), and
$A_{\lambda}$, $C_{\lambda}(\omega)$ and $S_{\lambda}(\omega)$ are
defined through (\ref{f1-s3}). We need to estimate $A_{\lambda}$
and $C_{\lambda}(\omega)$.

In this section we have so far computed $A_{\lambda}$ and
$C_{\lambda}, S_{\lambda}$ for $\lambda = 0$. We have from
(\ref{f7-s4}) and Lemma \ref{lem2-s4}
\begin{equation}\label{f9-s4}
A_0 = \frac{16}{15}, \ \ \ C_0(\omega) =  \frac{16 \sqrt{2} \pi
}{e^{-\frac{1}{2} \omega \pi} + e^{\frac{1}{2} \omega \pi}}.
\end{equation}
We caution that $\alpha_0 = \beta_0 = 1$ so for $\lambda = 0$ the
saddle point $(x, y) = (0, 0)$ is not dissipative. However, since
$-\alpha + \beta = - \lambda$ from (\ref{alpha-beta}), $(x, y) =
(0, 0)$ is a dissipative saddle for $\lambda >0$. We estimate
$A_{\lambda}$ and $C_{\lambda}(\omega)$ through $A_0,
C_0(\omega)$.

Let $U_{\varepsilon}$ be the $\varepsilon$-neighborhood of $(0,
0)$ in the $(x, y)$-plane, and $s^+_{\lambda}(\varepsilon)$ be the
time $\ell_{\lambda}$ enters $U_{\varepsilon}$. Similarly let
$-s^-_{\lambda}(\varepsilon)$ be the time $\ell_{\lambda}$ exits
$U_{\varepsilon}$. Then there exists $\varepsilon_0 > 0$ and $K$
independent of $\lambda$ such that for all $\lambda > 0$
sufficiently small, and all $s \in (s^+_{\lambda}(\varepsilon_0),
+\infty)$,
\begin{equation}\label{f11-s4}
|u_{\lambda}(s)| > 1- K \varepsilon_0, \ \ \ \ |v_{\lambda}(s)| <
K \varepsilon_0.
\end{equation}
Similarly we have, for all $s \in (-\infty,
-s^-_{\lambda}(\varepsilon_0))$,
\begin{equation}\label{f12-s4}
|v_{\lambda}(s)| > 1 - K\varepsilon_0, \ \ \ \ \ |u_{\lambda}(s)|
< K \varepsilon_0.
\end{equation}
Both (\ref{f11-s4}) and (\ref{f12-s4}) follow directly from the
fact that the local stable and the unstable manifold of equation
(\ref{f1-s1a}) are tangent to the $x$ and the $y$-axis
respectively at $(x, y) = (0, 0)$, and their curvatures are
bounded uniformly in $\lambda$. From (\ref{f11-s4}) and
(\ref{f12-s4}) we have, for all $s \in (-\infty,
-s^-_{\lambda}(\varepsilon_0)) \cup (s^+_{\lambda}(\varepsilon_0),
+\infty)$,
\begin{equation}\label{f13-s4}
\left| E_{\lambda}(s) \right| > \frac{1}{2}
\end{equation}
provided that $\lambda$ are sufficiently small. Let
$$
A_{\lambda, s^-, s^+} = \int_{-s^-}^{s^+}
(u_{\lambda}(s)+v_{\lambda}(s))(b_{\lambda}(s) +
a_{\lambda}(s))^2(b_{\lambda}(s)-a_{\lambda}(s)) e^{-\int_0^s
E_{\lambda}(t) dt} ds.
$$
We have from (\ref{f13-s4}) that
\begin{equation}\label{f14-s4}
|A_{\lambda} - A_{\lambda, s^-, s^+}| < \frac{1}{100}
\end{equation}
for all $\lambda$ that is sufficiently small provided that $s^- >
s^-_{\lambda}(\varepsilon_0), s^+ > s^+_{\lambda}(\varepsilon_0)$.

Let  $s^- = s^-_0(\frac{1}{2} \varepsilon_0)$, $s^+ =
s^+_0(\frac{1}{2} \varepsilon_0)$. From Proposition
\ref{prop1}(ii), we have $\ell_{\lambda}(-s^-),
\ell_{\lambda}(s^+) \in U_{\varepsilon_0}$ provided that $\lambda
\in (0, K^{-1}(\varepsilon_0))$ where $K(\varepsilon_0)>0$ is a
dependent of $\varepsilon_0$. We also have $\ell_{\lambda}(s) \cap
U_{\frac{1}{4} \varepsilon_0} = \emptyset$ for all $s \in (-s^-,
s^+)$, and this implies that for all $s \in (-s^-, s^+)$
\begin{equation}\label{f15-s4}
|u_{\lambda}(s) - u_0(s)|, \ \ \ |v_{\lambda}(s) - v_0(s)|<
K(\varepsilon_0) \lambda.
\end{equation}
(\ref{f15-s4}) and Proposition \ref{prop1}(ii) together implies
\begin{equation}\label{f16-s4}
|A_{\lambda, s^-, s^+} - A_{0, s^-, s^+}| < K(\varepsilon_0)
\lambda.
\end{equation}
That $A_{\lambda} > 1$ follows now from (\ref{f16-s4}),
(\ref{f14-s4}) and $A_0
> \frac{16}{15}$.

Estimates on $C_{\lambda}(\omega)$ is similar. Because
$C_{\lambda}(\omega)$ is a function in $\omega$, so would be the
constant $K(\varepsilon_0)$ in (\ref{f16-s4}). This is why
$\lambda_0$ in this proposition is a dependent of $R_{\omega}$.
\hfill $\square$

\end{document}